%% file: min.tex
\documentclass{amsart}

\usepackage{amsmath}
\usepackage{amssymb}
\usepackage{min}

\begin{document}

\input{mintop}

\input{min00} 
\input{min01} 
\input{min02} 
\input{min03} 
\input{min04} 
\input{min05} 
\input{min06} 
\bibliographystyle{plain}
\input{min.bbl}


\end{document}

%% file: mintop.tex
\title[Minimal Connected Simple Groups]
{Minimal Connected Simple Groups of finite Morley rank 
with Strongly Embedded Subgroups\\
{\vskip.1true in April 2005}}	

\author[J.~Burdges]
{%
Jeffrey Burdges\\
  Fakult\"at f\"ur Mathematik, Universit\"at Bielefeld
\email{burdges@math.rutgers.edu}%
}
\thanks{First author supported by NSF Graduate Research Fellowship,
and DFG at Wurzburg and Bielefeld, grant Te 242/3-1,
and acknowledges the
hospitality of the Universities of Birmingham and Manchester.}
\author[G.~Cherlin]{%
\\Gregory Cherlin\\
Department of Mathematics, Rutgers University
}
\thanks{Second author supported by NSF Grant DMS-0100794}
\author[E.~Jaligot]{%
\\Eric Jaligot\\
Equipe de Logique Math\'ematiques, Universit\'e de Paris VII
}
\thanks{All authors acknowedge the hospitality of the IGD at
Universit\'e Lyon I and support by the Newton Institute, Cambridge
during their Model Theory program, as well as the opportunity
afforded by the CIRM at Luminy for stimulating discussions.}
\maketitle
\thispagestyle{empty}

%% file: min00.tex

\section{Introduction}\label{sec:Intro}

The Algebraicity Conjecture states that a simple
group of finite Morley rank should be isomorphic with an algebraic group.
A program initiated by Borovik aims at controlling the 2-local structure in a
hypothetical minimal counterexample to the Algebraicity Conjecture.
There is now a large body of work on this program.
A fundamental division arises at the outset, according to
the structure of a Sylow 2-subgroup. In algebraic groups this structure
depends primarily on the characteristic of the base field. In groups of
finite Morley rank in general, in addition to the {\it even} and
{\it odd} type groups, which correspond naturally to the cases of
characteristic two or not two, respectively,  
we have two more cases, called {\it mixed} and
{\it degenerate} type. 
In the degenerate case the Sylow 2-subgroup is finite.

The cases of even and mixed type groups are well in hand, and it seems
that work in course of publication will show that the simple groups of
finite Morley rank of these two types are algebraic. Work on
degenerate type has hardly begun, though recently some interesting
approaches have emerged.  We deal here with odd type groups
exclusively.  In this context, the ``generic'' case is generally
considered to be that of groups of Pr\"ufer $2$-rank three or
more. The following result enables us to complete the analysis of the
generic case.

\begin{theorem}\label{main}
Let $G$ be a minimal connected simple group of finite Morley rank and of
odd type. Suppose that $G$ contains a proper definable strongly embedded
subgroup $M$. 
Then $G$ has Pr\"ufer $2$-rank one.
\end{theorem}

Combining this with known results, we will derive the following.

\begin{theorem}\label{application}
Let $G$ be a simple $K^*$-group of finite Morley rank of odd type, 
which is not algebraic. 
\begin{enumerate}
\item Then $G$ has Pr\"ufer $2$-rank at most two.
\item If $G$ is tame and minimal connected simple, 
and all the involutions in a standard Borel subgroup of $G$ are
central, then $G$ has Pr\"ufer $2$-rank one. 
\end{enumerate}

\end{theorem}

The proof of Theorem \ref{main} will be self-contained, while for Theorem
\ref{application} we will need to invoke an extensive body of material,
some in course of publication (or, alternatively, available in
\cite{Bu-Th}). We will now enlarge on the terminology used above.

The Pr\"ufer $2$-rank of a group of finite Morley
rank is the Pr\"ufer rank of a maximal 2-divisible abelian subgroup,
and is always finite.

A group of finite Morley rank is {\it tame} if it involves no bad
field, and is a $K^*$-group if every proper definable infinite simple
section is algebraic.  The Borovik program was initially directed
primarily at tame $K^*$-groups.  It is known that a {\it tame}
$K^*$-group of finite Morley rank and odd type has Pr\"ufer $2$-rank
at most two, which gives a qualitative version of the Borovik program
in the tame case; most of the analysis does not depend on tameness,
up to the point of a reduction to the minimal connected simple case, 
treated in \cite{CJ01} by methods specific to the tame case.
Each clause of Theorem \ref{application} improves
on this result; the first clause eliminates the tameness hypothesis,
while the second clause eliminates one of the configurations in
Pr\"ufer rank two which was left open in \cite{CJ01}.

The first clause of Theorem \ref{application} completes the analysis of the
``generic'' $K^*$-group of finite Morley rank of odd type. An outline of
the various special cases which require further analysis, and some ideas as
to how they may be approached, is given in Chapter 11 of \cite{Bu-Th}.

In the next section we will give the derivation of Theorem
\ref{application} from Theorem \ref{main}, and we will also make some
preliminary remarks concerning the proof of Theorem \ref{main}.
Subsequent sections are devoted to the proof of Theorem \ref{main},
which is divided into two cases. These two cases will be treated by
quite distinct methods. The first is handled relativel rapidly in \S4,
while the other occupies us for another three sections.  Our method of
analysis is heavily influenced by \cite{CJ01}, which obtained similar
results under the hypothesis of tameness, and also gave an analysis of
the problematic configurations in Pr\"ufer $2$-ranks one and two (with
two such configurations in each rank). As noted, our theorem also
eliminates one of these four problematic configurations.

A case division exploited throughout \cite{CJ01}, and which we
continue to make use of, is the following.
Let a {\it standard} Borel subgroup be a Borel subgroup which contains a
Sylow\oo\ 2-subgroup of $G$. Then our case division is as follows: either
(I) all the involutions of a standard Borel subgroup are central, or
(II) not.
In the former case, the method used in \cite{CJ01} in high Pr\"ufer $2$-rank is
inappropriate outside the tame case, and is replaced here by another method
which, as noted, yields a sharper result in Pr\"ufer $2$-rank two. In the
latter case, one becomes involved in a close analysis of intersections of
Borel subgroups. The basic lemma concerning such intersections 
in the tame case is the following.

\begin{fact}[{\cite[\S3.4]{CJ01}}]\label{jaligotlemma}
Let $G$ be a tame minimal connected simple group of odd type and finite
Morley rank. Assume that $B_1$ and $B_2$ are distinct Borel subgroups of
$G$ with $O(B_1),O(B_2)\ne 1$. Then $F(B_1)\intersect F(B_2)=1$. 
\end{fact}

An important consequence of this lemma is that $B_1\intersect B_2$ is
abelian under the stated hypotheses (as $B_1/F(B_1)$ and $B_2/F(B_2)$
are abelian). It would be equally useful to have $(B_1\intersect
B_2)\oo$ abelian.  We see no direct way of proving anything as strong
as this in general, but we continue to work in this general direction,
at the price of a much more elaborate analysis.  This was initiated in
\cite{Bu-Th}, whose Theorem 9.2 gives a very elaborate analog of
Fact \ref{jaligotlemma}, describing the configuration that arises from
the situation in which $(B_1\intersect B_2)\oo$ is nonabelian.  

In the tame case, in high Pr\"ufer $2$-rank, the case in which a
standard Borel subgroup $B$ contains an involution not in its center
was disposed of in short order via the powerful Fact
\ref{jaligotlemma}, whereas the analysis of the other case was quite
long. Here it appears the situation is entirely reversed, with the
case of central involutions being dispatched relatively expeditiously
via the construction of two disjoint generic subsets.  However, in
\cite{CJ01} a portion of the preliminary analysis of the case of
central involutions actually went toward showing that the hypotheses
of our main theorem are satisfied in that case.


%% file: min01.tex

\section{Preliminaries}\label{preliminaries}

\subsection{Theorem \protect\ref{application}}
\label{subsec:application}

\medskip
We will now discuss the prior results which link Theorems \ref{main} and
\ref{application}.
The two facts which require explanation are as follows.

\begin{fact}\label{genericoddtype}
Let $G$ be a simple nonalgebraic 
$K^*$-group of finite Morley rank and odd type,
and Pr\"ufer $2$-rank at least three.
Then $G$ is minimal connected simple, and has a proper definable
strongly embedded subgroup. 
\end{fact}

\begin{fact}\label{tameoddtype}
Let $G$ be a tame and minimal connected simple group of finite Morley
rank and odd type, and Pr\"ufer $2$-rank at least two,
and suppose that there is a standard Borel subgroup $B$ 
of $G$
such that every involution of $B$ lies in $Z(B)$.
Then $G$ has a proper definable strongly embedded 
subgroup. 
\end{fact}

We deal first with Fact \ref{genericoddtype}. 

A general reference for this fact, and for a great deal of prior material,
is the thesis \cite{Bu-Th}, notably Theorem 10.15 in that reference,
combined with Theorems 8.14 and 8.18. The latter two theorems
may be found in \cite{BBN-UC}, while Theorem 10.15 will appear in
\cite{Bu-BC}.
These results depend in turn on a very substantial body of material due to
Borovik and others. For the reader's convenience we give some
additional details.

\begin{namedtheorem}{Theorem 10.15}
Let $G$ be a simple $K^*$-group of finite Morley rank and odd type with
Pr\"ufer 2-rank at least three. Then either $G$ has a proper 2-generated
core, or $G$ is an algebraic group over an algebraically closed field
of characteristic not two.
\end{namedtheorem}

The {\em 2-generated core} is
defined as follows. Let $S$ be a Sylow 2-subgroup of $G$. Let
$\Gamma_{2,S}^e(G)$ be the definable closure of the subgroup of $G$ generated
by all subgroups of the form 
$$\hbox{$N_G(V)$ (for $V\le S$, $V$ an elementary abelian 2-group of
  rank two)}$$ 
Then $\Gamma_{2,S}^e(G)$ is the {\it 2-generated core} of $G$ with respect to
$S$, and as $S$ varies the groups $\Gamma_{2,S}^e$ are conjugate. So the
2-generated core of $G$ is well-defined up to conjugacy. 

Having a proper 2-generated core is weaker than having a strongly embedded
subgroup, but is the first step in this direction. The next step is given
by the following, found in \cite{BBN-UC}; see also
\cite[Chap.~8]{Bu-Th}, 

\begin{namedtheorem}{Strong Embedding Theorem}
Let $G$ be a simple $K^*$-group of finite Morley rank and odd type, and
suppose that $G$ has normal $2$-rank at least three and Pr\"ufer $2$-rank
at least two. Let $S$ be a Sylow 2-subgroup of $G$, and $M=\Gamma_{2,S}^e(G)$
the associated 2-generated core. If $M$ is a proper subgroup of $G$, then
$G$ is a minimal connected simple group and $M$ is strongly embedded in
$G$. Furthermore $M\oo$ is a Borel subgroup of $G$ and $S$ is connected.
\end{namedtheorem}

Note that Fact \ref{genericoddtype} now follows.

We turn next to Fact \ref{tameoddtype}, where we deal with the tame case.
This fact concerns one of the four configurations which were
considered explicitly in \cite{CJ01}, in the 
context of tame minimal connected simple groups of finite Morley rank and
odd type, and which were not eliminated there. It shows that the
strong embedding hypothesis with which we will be working is valid in
that case.
The groups under consideration in \cite{CJ01} are tame groups of
finite Morley rank, minimal among connected simple groups.
The configuration in which all involutions in a standard Borel subgroup
$B$ are central, and with the additional hypothesis of Pr\"ufer rank
at least two, is taken up in the first part of \S7 of \cite{CJ01}.

There it is shown, first, that a standard Borel subgroup $B$ is nilpotent, and
then, after some further analysis, that the normalizer of $B$ is
strongly embedded.  It is also clear in this case that each Sylow\oo\
2-subgroup is contained in a unique standard Borel subgroup, and in
particular the standard Borel subgroups are conjugate, so they all
have the same properties.

Note that the case division considered in \cite[\S\S6,7]{CJ01} is
equivalent to the following: a standard Borel subgroup is, or is not,
strongly embedded. In all cases that survive analysis the
centralizer of an involution contains a standard Borel subgroup. 
Here we will begin with a strongly embedded subgroup, and
then see quickly that the normalizer of 
a standard Borel subgroup is strongly
embedded; for us, the interesting case division lies farther on, 
within this configuration, and the harder of the two cases is one that
vanished after a few lines of analysis in the tame context.

\subsection{The context of Theorem \protect \ref{main}}

We now turn to Theorem \ref{main}. 
Suppose that $G$ is a minimal connected simple group
of finite Morley rank and odd type, and that $M$ is a proper definable
strongly embedded subgroup of $G$. 
We also assume that 
$$\hbox{$G$ has Pr\"ufer $2$-rank at least two.}\leqno(\dag)$$
We set $B=M\oo$, and to justify this notation we prove the following.

\begin{lemma}
$B$ is a standard Borel subgroup of $G$, and $M=N(B)$.
\end{lemma}
\begin{proof}

Since $M$ is a strongly embedded subgroup of $G$, it contains a Sylow
2-subgroup of $G$, and hence $B$ contains a Sylow\oo\ 2-subgroup $S$
of $G$. By the minimality of $G$, $B$ is solvable. 

Now suppose $B\le H$ with $H$ a Borel subgroup of $G$. We claim that
$B=H$. Let $V\le\Omega_1(S)$ have order four. 
As $V\le H$, in particular $V$ acts on $H$, and as $H$ is connected
solvable of odd type, a fundamental generation property given
in \cite[5.14]{Bo95} implies that 
$$H=\<C_H\oo(v):v\in V^\#\>$$
On the other hand by strong embedding of $M$, we have
$C_H\oo(v)\le B$ for all involutions $v\in V^\#$, and thus $H=B$.
So $B$ is a standard Borel subgroup of $G$.

As $B=M\oo$, we have $M\le N(B)$. Conversely, 
as $B\le M$ contains involutions and $M$ is
strongly embedded, we conclude $N(B)\le M$. Thus $M=N(B)$.
\end{proof}

Part of the foregoing argument can be strengthened as follows
to a ``black hole'' principle (the term goes back to Harada).
We record this for future reference.

\begin{lemma}\label{blackhole}
Let $H$ be a connected definable proper subgroup of $G$ and $V$ an
elementary abelian $2$-subgroup of $B$ of rank 2. 
If $V$ normalizes $H$ then $H\le B$.
\end{lemma}

Now we intend to make a case division based on whether or not all the
involutions of $B$ are central in $B$. As $M$ is strongly embedded in
$G$, all involutions of $M$ are conjugate in $M$. Hence all
involutions of $M$ lie in $B$, and if one involution is central in $B$
then all involutions of $B$ lie in its center.
So our case division is actually the following: 
either all involutions of $B$ are central in $B$, or none are. 

Before entering into the analysis of individual cases,
we give a genericity result that holds in both cases.
This was given already in \cite{Bu-Th} and somewhat more explicitly in 
\cite{BBN-UC}, under stronger hypotheses which are not actually used
at this point in the argument. For the reader's convenience we indicate 
the gist of the argument, which depends on the following generic
covering lemma, a result for which we will have further use in
\S\ref{sec:invariance}. 

\begin{fact}[{\cite[3.3]{CJ01}}]\label{genericcovering}
Let $G$ be a connected group of finite Morley rank, and $B$ a definable
subgroup of finite index in its normalizer. Suppose that there is a
non-generic subset $X$ of $B$ such that $B\intersect B^g\includedin X$ for
$g\in G\setminus N(B)$. Then $\Union_{g\in G}B^g$ is generic in $G$.
\end{fact}

\begin{lemma}\label{Bgenericcovering}
The union $\Union_{g\in G}B^g$ is a generic subset of $G$.
\end{lemma}

\begin{proof}
By strong embedding, $M\intersect M^g$ contains no involutions for
$g\notin M$, and in particular $B\intersect B^g$ is a $2^\perp$-group
for $g\notin M$. 
By a very general lemma given as
\cite[8.17]{Bu-Th} or \cite[5.10]{BBN-UC}, 
the union $\Union_{g\in   G\setminus N(B)}({B\intersect B^g})$ is 
contained in a proper definable 
subgroup of $B$, in fact a $2^\perp$-group,
 and hence is not generic in $B$. 
Then by the generic covering lemma 
the claim follows.
\end{proof}

In the remainder of the article, we will prove Theorem \ref{main},
dividing the analysis into two very different cases, which rely on
very different methods. We first summarize some very general
background material, used throughout, and sufficient for the treatment
of the first of our two cases. The second case will involve a further 
body of material which is both more extensive, and of more recent vintage.


%% file: min02.tex

\section{Background material}
\label{background}

The material of the present section is for the most part well
known. The most subtle item comes from the theory of Carter subgroups
and is due to Fr\'econ: Fact \ref{carter-hall}, and its consequence
Lemma \ref{carter-sylow}, below.

\subsection{Unipotence}

\begin{definition}
Let $p$ be a prime, $G$ a group of finite Morley rank, and $P$ a subgroup
of $G$. 
\begin{enumerate}
\item $P$ is said to be {\em $p$-unipotent} if $P$ is a solvable connected
definable $p$-subgroup of $G$ of finite exponent.
\item $U_p(G)$ is the largest normal $p$-unipotent subgroup of $G$.
\end{enumerate}
\end{definition}

The group $U_p(G)$ is well-defined, by an elementary
argument.

\begin{fact}[{\cite{Ne90}}]\label{Up}
Let $H$ be a solvable group of finite Morley rank and $P$ a $p$-unipotent
subgroup of $H$. Then $P\le U_p(H)$.
\end{fact}

This is phrased somewhat differently in \cite{Ne90}. The essential
point is that $H\oo/F\oo(H\oo)$ is divisible abelian. From this it
follows that $P\le F\oo(H)$, and then the structure of nilpotent
groups applies, as in \cite[\S6.4]{BN}.

The next lemma is a weak form of Fact \ref{jaligotlemma}, and has a
similar proof, given also in \cite{Bu-GCS}.  The virtue of this lemma
is that it holds without any assumption of tameness. While this will
suffice for the purposes of the next section, we will need much more
subtle variations subsequently.

\begin{lemma}\label{punipotent}
Let $G$ be a minimal connected simple group of finite Morley rank, $p$ a
prime, and $P$ a nontrivial $p$-unipotent subgroup of $G$.
Then $P$ is contained in a unique Borel subgroup of $G$.
\end{lemma}
\begin{proof}
Suppose on the contrary that $B_1$ and $B_2$ are two Borel subgroups
containing $P$, and chosen so that $Q=U_p(B_1\intersect B_2)$ is maximal.
Then $Q\le U_p(B_1)$. 

Suppose $Q<U_p(B_1)$. Then $N_{U_p(B_1)}\oo(Q)>Q$
by the normalizer condition, \cite[\S6.4]{BN}. Put $N\oo(Q)$ in a Borel
subgroup $B_3$; by the maximality of $Q$, we find $B_3=B_1$. Then $B_3\ne
B_2$, so we must have $Q=U_p(B_2)$. Then $B_2\le N\oo(Q)\le B_1$, a
contradiction.

There remains the possibility that $Q=U_p(B_1)=U_p(B_2)$. But then
$B_1,B_2=N\oo(Q)$, again a contradiction.
\end{proof}

\subsection{Sylow subgroups}

\begin{lemma}\label{psylow}
 Let $P$ be a Sylow $p$-subgroup of a connected 
solvable group of
 finite Morley rank. Then $P$ decomposes as $U*T$, a
 central product, with $U$ $p$-unipotent 
 and $T$ a $p$-torus (that is, a divisible abelian
 $p$-group).
\end{lemma}
\begin{proof}
This is essentially \cite[\S6.4]{BN}, where a similar structure
theorem is given more generally for $p$-subgroups of solvable groups
of finite Morley rank. When $H$ is in addition connected, 
then by \cite[9.39]{BN}, 
its Sylow (or indeed its Hall) subgroups are connected, and this
gives the connectivity of the factor $U$.
\end{proof}

\begin{lemma}\label{sylow:basics}
Let $G$ be a group of finite Morley rank, $P$ a Sylow $2$-subgroup, 
$H$ a normal subgroup.
Then:
\begin{enumerate}
\item $P\intersect H$ is a Sylow $2$-subgroup of $H$.
\item If $H$ is definable and $\bar G=G/H$, 
then $\bar P$ is a Sylow $2$-subgroup of $\bar G$.
\end{enumerate}
\end{lemma}

The first follows directly from the conjugacy of Sylow $2$-subgroups.
The second point is given in \cite{PoWa-LS}.

\subsection{Genericity}

\begin{fact}[{\cite[3.4]{CJ01}}]\label{N(B)}
Let $G$ be a connected group of finite Morley rank. Suppose that $B$
is a definable subgroup of finite index in its normalizer such that
$\Union_{g\in G} B^g$ is generic in $G$. Suppose that $x\in
N_G(B)\setminus B$.  Let $X$ be the set
$$\{x'\in xB:\hbox{$x'\in (\<x\>B)^g$ for some $g\in G\setminus N(B)$}\}$$
Then $X$ is generic in $xB$.
\end{fact}

\begin{fact}[{\cite[3.6]{CJ01}}]\label{genericorder}
Let $H$ be a group of finite Morley rank such that $H\oo$ is abelian, and
let $xH\oo$ be a coset whose elements are generically of fixed order $n$.
Then every element of $xH\oo$ is of order $n$.
\end{fact}

\begin{lemma}\label{punipotentaction}
 Let $H$ be a group of finite Morley rank with $H\oo$
 solvable, and $H/H\oo$ of prime order $p$. Suppose the
 elements of every coset of $H\oo$ other than $H\oo$ are
 generically of order $p$. If some element of $H\setminus
 H\oo$ has an infinite centralizer in $H\oo$, then $H\oo$
 contains a nontrivial $p$-unipotent subgroup.
\end{lemma}
\begin{proof}
This was proved under the assumption that $H\oo$ is
nilpotent in \cite[3.7]{CJ01}.
However, under the stated hypotheses (or slightly weaker
ones; one such coset suffices) the solvability of $H\oo$
implies its nilpotence, by \cite[Cor.~16]{JW00}.
\end{proof}

\subsection{Carter subgroups}

A Carter subgroup of a solvable group of finite Morley rank
is a definable self-normalizing nilpotent subgroup. 
If $H$ is a connected solvable group of finite Morley rank,
then by \cite[5.5.10, 5.5.12]{Wa-NCC} and 
\cite[1.1]{Fre00},
it has a Carter subgroup, and any two such are conjugate;
and by \cite[3.2]{Fre00},
its Carter subgroups are connected. 

\begin{fact}[{\cite[3.2]{Fre00}}]\label{nearcarter}
Let $H$ be a connected solvable group of finite Morley rank, and $Q$ a
nilpotent subgroup of $H$ with $[N_H(Q):Q]$ finite. Then $Q$
is a Carter subgroup of $H$.
\end{fact}

This includes the more elementary result
that the centralizer of any element of finite order in a connected
solvable group is infinite (compare \cite{Ja-FFG}).

\begin{fact}[{\cite[7.15]{Fre00}}]\label{carter-hall}
Let $H$ be a connected solvable group of finite Morley rank, and $R$ a Hall
$\pi$-subgroup of $H$ for some set $\pi$ of primes. Then $N_H(R)$ contains
a Carter subgroup of $H$.
\end{fact}

This has the following important consequence.

\begin{lemma}\label{carter-sylow}
Let $H$ be a connected solvable group of finite Morley rank with
$U_p(H)=1$, and let $Q$ be a Carter subgroup of $H$. Then $Q$ contains a
Sylow $p$-subgroup of $H$.
\end{lemma}
\begin{proof}
Let $S$ be a Sylow $p$-subgroup of $H$. 
By Fact \ref{carter-hall}, $N_H(S)$
contains a Carter subgroup $Q_0$ of $H$. $Q_0$ is connected.
It follows easily from Lemma \ref{psylow} that $S$ is a $p$-torus, that is,
abelian and $p$-divisible. So $Q_0\le N\oo(S)=C\oo(S)$.
Thus $S\le N(Q_0)=Q_0$. Now as $Q_0$ and $Q$ are conjugate, our claim
follows.
\end{proof}


%% file: min03.tex

\section{Central involutions}
\label{sec:se}

\subsection{The setup}
We take up the proof of Theorem \ref{main}. We dispose of one case in
the present section, and the other will occupy us to the end of the
paper.

So $G$ is a minimal connected simple group of finite
Morley rank and odd type, with $M$ definable and strongly embedded,
and we assume that 
$$\hbox{$G$ has Pr\"ufer $2$-rank at least two.}\leqno(\dag)$$
As explained in \S\ref{preliminaries}, we then have
$M=N(B)$ strongly embedded, with $B$ a standard Borel subgroup.

In the present section we take up the first of our two cases, namely:
$$\hbox{All involutions in $B$ are central in $B$}\leqno(\hbox{\bf Case I})$$
In particular $C\oo(i)=B$ for each involution $i$ in $B$.

We will derive a contradiction in this case by constructing
two disjoint generic subsets of $G$. That is, in Case I the Pr\"ufer
rank of $G$ is at most $1$.

The two generic subsets in question will be $BI_1$ and $BC(\sigma)B$
where $I_1={I(G)\setminus I(M)}$ and $\sigma\in M\setminus B$.
We must deal with the following issues: that the ranks
of the two sets in question are the same as the ranks of the
corresponding Cartesian products $B\times I_1$ and $B\times
C(\sigma)\times B$; that these ranks coincide with the rank of $G$;
and that the sets in question are disjoint. We must also produce a
suitable element $\sigma$, but that is immediate since the involutions
of $B$ form an elementary abelian subgroup which by hypothesis has 2-rank
at least two, and they are conjugate under the action of $M$, by
strong embedding, and are central in $B$.

That the set $BI$ or $BI_1$ turns out to be generic in $G$ is
perfectly natural, but in the case of the set $BC(\sigma)B$ it is
surprising.

We record this notation.

\begin{notation}
\

1. $I=I(G)$ and $I_1=I\setminus M=I\setminus B$.

2. $\sigma\in M\setminus B$ (fixed).
\end{notation}

\subsection{The first generic subset}

Recall that an element $a$ of $G^\#$ is said to be {\it strongly real} if
it is a product of two involutions, in which case $a$ is inverted by
these two involutions. 

The following fundamental fact will be used in
both of our cases.

\begin{fact}\cite[10.19]{BN}\label{BN:stronglyreal}
In a group of finite Morley rank with a
definable strongly embedded subgroup $M$, if $a$ is a strongly real
element commuting with an involution in $M$, then every involution
which inverts $a$ lies in $M$. 
\end{fact}

In Case I,
this applies to every strongly real element of $B$,
so we arrive at the following under our present assumptions.

\begin{lemma}\label{stronglyreal}
The strongly real elements of $B$ are its involutions.
\end{lemma}
\begin{proof}
By Fact \ref{BN:stronglyreal}, if $b\in B$ is strongly real and $j$ is an
involution inverting $b$, then we have $j\in B$. 
By our case
hypothesis $j$ commutes with $b$, so
$b$ is an involution.

Conversely, any involution in $B$ lies in $Z(B)$ and is the product of
two involutions by the hypothesis $(\dag)$.
\end{proof}

\begin{lemma}\label{1stgeneric}
The set $BI_1$ is generic in $G$, and the multiplication map 
$$B\times I_1\to G$$ 
is injective.
\end{lemma}
\begin{proof}
Notice first that $I(B)$ is finite since $M/B$ operates transitively
on this set. Hence $\rk(I)=\rk(I_1)$.

Let $i\in I(B)$. Recall that $C\oo(i)=B$, and $\rk(I)=\rk(G/C(i))$ as
all involutions of $G$ are conjugate. 
So we have $\rk(G)=\rk(C(i))+\rk(G/C(i))=\rk(B)+\rk(I)=\rk(B)+\rk(I_1)$.
It suffices therefore to check that the multiplication map
$B\times I_1\to G$ is injective.

Supposing the contrary, we have a nontrivial intersection
$bI_1\intersect I_1$ with $b\in B^\#$, which yields an equation $b=jk$
with $j,k\in I_1$.  Thus $b\in B$ is strongly real.  By Lemma
\ref{stronglyreal}, $b$ is an involution, and $j$ centralizes $b$,
hence lies in $M$, a contradiction.
\end{proof}

From this genericity result, others of the same type can be deduced,
for sets of the following form.

\begin{notation}
For $g\in G$, we set $I_g=(B g I)\intersect I_1$
\end{notation}

\begin{lemma}\label{2ndgeneric}
For any $g\in G$, the set $B g I$ is
generic in $G$, and the set $I_g$ is generic in $I$.
\end{lemma}
\begin{proof}
The set $BI$ is generic in $G$. Conjugating by $g$, the set $B^gI$
is generic in $G$. Translating on the left by $g$, the set $B g I$ is
generic in $G$. Hence also $(B g I)\intersect (BI_1)$ is generic in
$G$.

Now $BI_g=(B g I)\intersect (BI_1)$ is
generic in $G$. Hence $\rk(G)=\rk(BI_g)=\rk(B)+\rk(I_g)$ as the
multiplication map restricted to $B\times I_g$ is injective (Lemma
\ref{1stgeneric}), and thus 
$\rk(B)+\rk(I_g)=\rk(B)+\rk(I)$ and $\rk(I_g)=\rk(I)$.
\end{proof}

\subsection{Disjointness and centralizers}

The rest of the argument requires some more structural information. 

\begin{lemma}
The intersection of two distinct conjugates of $B$ is finite.
\end{lemma}
\begin{proof}
Let $S$ be a Sylow 2-subgroup of $B$, and $V=\Omega_1(S)$. Then $V\le
Z(B)$. Hence $V$ centralizes $B\intersect B^g$.

Suppose the group $(B\intersect B^g)\oo$ is nontrivial.  Then it has a
nontrivial Carter subgroup $Q$.  As $V$ centralizes $Q$, it normalizes
$N\oo(Q)$. By Lemma \ref{blackhole}, $N\oo(Q)\le B$. Hence
$N_{B^g}\oo(Q)\le (B\intersect B^g)\oo$ and thus $N_{B^g}\oo(Q)=Q$.
By Fact \ref{nearcarter}, $Q$ is a Carter subgroup of $B^g$. Then by
Lemma \ref{carter-sylow}, $Q$ contains a Sylow 2-subgroup of $B^g$.
In particular $B\intersect B^g$ contains an involution $i$ and hence
$B=C\oo(i)=B^g$.
\end{proof}

\begin{lemma}
The intersection of any two distinct conjugates of $B$ is trivial.
\end{lemma}
\begin{proof}
Supposing the contrary, let $B\intersect B^g$ be finite and
nontrivial, and let
$x\in B\intersect B^g$ have prime order
$p$. Let $P$ be a Sylow $p$-subgroup of $B$ containing $x$.  By
Lemma \ref{psylow}, $P$ is a central product
of the form $U*T$ with $U$ connected, nilpotent of bounded exponent
and $T$ divisible abelian.  Here $d(U)$ is also a connected nilpotent
$p$-group centralizing $T$, so $P=d(U)*T$ and $U=d(U)$ is definable.

If $U\ne 1$, let $U_0=C_U\oo(x)$. Then $U_0\ne1$ (\cite[6.20]{BN}).
Now $U_0$ is contained in a unique Borel subgroup by Lemma
\ref{punipotent}, and hence $C\oo(x)\le B$. But as $B$ and $B^g$ are
conjugate, a Sylow $p$-subgroup of $B^g$ has the same form, and hence
this argument shows $C\oo(x)\le B^g$ as well. But then $U_0\le B^g$
and $B=B^g$.

This contradiction shows that $U=1$ and hence $P$ is a $p$-torus.
Then Lemma \ref{carter-sylow} shows that there is a Carter subgroup
$Q$ of $B$ containing $P$. Furthermore, the same lemma shows that $Q$
contains a Sylow 2-subgroup $S$ of $B$. As $Q$ is nilpotent, it
follows that $x$ (in $P$) commutes with $S$. Similarly, $x$ commutes
with a Sylow 2-subgroup of $B^g$.

Let $V$ be a four-group contained in $S$. Then $V$ normalizes
$C\oo(x)$ and hence by Lemma \ref{blackhole} we have $C\oo(x)\le B$,
and in particular $B$ contains a Sylow 2-subgroup of $B^g$, forcing
$B=B^g$, a contradiction.
\end{proof}

\begin{lemma}\label{centralizerfinite}
Let $x\in N(B)\setminus B$. Then the centralizer $C_B(x)$ is finite.
\end{lemma}
\begin{proof}
Suppose $C_B(x)$ is infinite. Replacing $x$ by a power, we may
suppose that the order $p$ of $x$ modulo $B$ is a prime.

Consider the subset $X$ of $xB$ defined as 
$$\{x'\in xB:\hbox{$x'\in (\<x\>B)^g$ 
for some $g\in G\setminus N(B)$}\}$$
This is generic in
$xB$ by Fact \ref{N(B)}.  For $x'\in X$ we have ${x'}^p\in B\intersect
B^g=1$ for some $g\in G\setminus N(B)$. The same applies to any coset
in $\<x\>B$ other than $B$.

By Lemma \ref{punipotentaction} and Fact \ref{Up}, $U_p(B)>1$. Then for $x'\in
X$, with $x'\in N(B^g)$, $g\in G\setminus N(B)$, consider the action
of $x'$ on $U_p(B)$ and $U_p(B^g)$. It follows \cite[6.20]{BN} that
the groups $U_p(C_B(x'))$ and $U_p(C_{B^g}(x'))$ are nontrivial. Letting
$B_1$ be a Borel subgroup containing $C\oo(x')$, 
it follows by Lemma \ref{punipotent} that $B=B_1=B^g$, a contradiction.
\end{proof}

This yields strong information concerning the coset $B\sigma$
(which is in fact an arbitrary coset of $B$ in $M$, other than $B$).

\begin{lemma}\label{sigmastronglyreal}
\ 
\begin{enumerate}
\item $B\sigma=\sigma^B$.
\item The elements of $B\sigma$ are all strongly real, and lie outside
$\Union_{g\in G}B^g$.
\end{enumerate}
\end{lemma}
\begin{proof}
\ 

1. Since $\sigma\in N(B)$ we have $\sigma^B\includedin B\sigma$.
Furthermore, as $C_B(\sigma)$ is finite, we have
$\rk(\sigma^B)=\rk(B)=\rk(B\sigma)$.  As $B\sigma$ has Morley degree
one, it consists of a single $B$-conjugacy class.

2. As the set $I_\sigma$ is generic in $I$ by Lemma \ref{2ndgeneric}, it is
nonempty, and we have an equation $xi=j$ with $x\in B\sigma$ and
$i,j\in I$. Hence $x$ is strongly real, and since $B\sigma$ is a
single $B$-conjugacy class, all of its elements are strongly real.

Now the strongly real elements of $B$ are involutions. Hence no
element of $B\sigma$ can be conjugate to an element of $B$.
\end{proof}

\begin{lemma}\label{inversion}
Let $g\in G$ be strongly real, inverted by the involution $i$, with
$g$ not an involution. Then $i$ acts on $C\oo(g)$ by inversion.
\end{lemma}
\begin{proof}
We may suppose $i\in B$. It suffices to show that $C\oo(i,g)=1$
\cite[p.~79, Ex.~13,15]{BN}.
We have $C\oo(i,g)\le C\oo(i)=B$, and similarly $C\oo(i,g)\le B^g$. 
So if $C\oo(i,g)\ne 1$ then
$g\in N(B)$, and we contradict Lemma \ref{centralizerfinite} or Lemma
\ref{stronglyreal}. 
\end{proof}

\subsection{The second generic set}

We now take up the proof that $BC(\sigma)B$ is generic in $G$
and disjoint from $BI_1$. Let us begin with the latter point

\begin{lemma}\label{M/B}
$M/B$ has odd order.
\end{lemma}
\begin{proof}
If not, we may choose $\sigma$ of order $2$ modulo $B$.
By Lemma \ref{sylow:basics} we may choose $\sigma$ to be a
$2$-element. Let $i$ be an involution in the cyclic group $\<\sigma\>$.
Then $i\in B$. 

By Lemma \ref{sigmastronglyreal} the element $\sigma$ is
strongly real. Let $j$ be an involution inverting $\sigma$.
As $j\in C(i)$, we have $j\in M$, so $j\in B$.
But $j\sigma$ is another involution in $M$, 
so $j\sigma\in B$. Hence
$\sigma\in B$, a contradiction.
\end{proof}

\begin{lemma}\label{C(sigma),I1}
$C(\sigma)\intersect I_1=\emptyset$
\end{lemma}
\begin{proof}
In view of Lemma \ref{sigmastronglyreal},
there is an involution $i$ inverting $\sigma$.
Suppose toward a contradiction that there is also some
$j\in C(\sigma)\intersect I_1$, and let $M^g$ be the conjugate of
$M$ containing $j$. As $\sigma$ is strongly real, Fact
\ref{BN:stronglyreal} implies that $i\in M^g$, hence $i\in B^g$.
But $\sigma$ normalizes
$B^g=C\oo(j)$, so $\sigma\in M^g$. Now $\sigma^2=[i,\sigma]\in B^g$,
so $\sigma\in B^g$ by Lemma \ref{M/B}, and this contradicts Lemma
\ref{sigmastronglyreal} (2).

\end{proof}

\begin{lemma}\label{disjoint}
$BC(\sigma)B$ and $BI_1$ are disjoint.
\end{lemma}
\begin{proof}
Supposing the contrary, we have $b_1cb_2\in I_1$ for some $b_1,b_2\in B$
and $c\in C(\sigma)$, and conjugating by $b_1$ gives 
$$cb\in I_1\leqno(1)$$
with $b\in B$. Conjugating by $\sigma$ gives
$$cb^\sigma\in I_1\leqno(2)$$
It follows that $b^{-1}b^\sigma=(cb)^{-1}(cb^\sigma)$ 
is inverted by an element $j$ of $I_1$. 

If $b\in C(\sigma)$, then $cb\in C(\sigma)\intersect I_1$, 
contradicting Lemma \ref{C(sigma),I1}.
So $b^{-1}b^\sigma$ is nontrivial, and is a strongly real element of $B$.
So $j\in M$ by Fact \ref{BN:stronglyreal},
hence $j\in B$, which is a contradiction since $j\in I_1$.
\end{proof}

Now we compute the rank of $BC(\sigma)B$, getting the expected value.

\begin{lemma}\label{rankBCB}
The set $BC(\sigma)B$ has rank $2\rk(B)+\rk(C(\sigma))$.
\end{lemma}
\begin{proof}
Let $C_\sigma$ be $C(\sigma)\setminus N(B)$.  As $C_B(\sigma)$ is
finite by Lemma \ref{centralizerfinite}, 
the set $C_\sigma$ differs from $C(\sigma)$ by a finite set.
On the other hand, the centralizer $C_G(\sigma)$ is infinite, as
otherwise the conjugacy class $\sigma^G$ would be generic in $G$ and
hence meet $\Union_{g\in G}B^g$, contradicting Lemma
\ref{sigmastronglyreal}. 
So $\rk(C_\sigma)=\rk(C(\sigma))$.

It now suffices to check that the multiplication map
$$\mu:B\times C_\sigma\times B\to G$$
has finite fibers.

If $g=b c b'$ with $b,b'\in B$, $c\in C_\sigma$, and $\mu^{-1}(g)$ is
infinite, then the same applies to $c=b^{-1}g{b'}^{-1}$.
So we consider $\mu^{-1}(c)$ with $c\in C_\sigma$ fixed. That is, we
examine the solutions to the equation
$$b_1c'b_2=c\leqno(*)$$
with $b_1,b_2\in B$ and $c'\in C_\sigma$.
Applying $\sigma$, we get
$$b_1^\sigma c' b_2^\sigma=c$$
Combining these two equations yields
$$c'=b_1'c'b_2'$$
with $b_1'=b_1^{-1}b_1^\sigma$, $b_2'=b_2^\sigma b_2^{-1}$.
So $(b_1')^{c'}={b_2'}^{-1}\in B\intersect B^{c'}$, and as $c'\notin
N(B)$ we have 
$b_1'=b_2'=1$, and $b_1,b_2\in C_B(\sigma)$, which allows finitely many
possibilities, by Lemma \ref{centralizerfinite}.
\end{proof}

It remains now to recompute the rank of $G$. 
We have already noticed that $\rk(G)=\rk(B)+\rk(I)$, but we now need
$\rk(G)=2\rk(B)+\rk(C(\sigma))$, or in other words 
$$\rk(I)=\rk(B)+\rk(C(\sigma))$$ 
We know $\rk(I_\sigma)=\rk(I)$, so it will suffice to show that
the rank of $I_\sigma$ is $\rk(B)+\rk(C(\sigma))$.

\begin{lemma}\label{rankG}
$\rk(G)=2\rk(B)+\rk(C_G(\sigma))$
\end{lemma}
\begin{proof}
We have $\rk(G)=\rk(B)+\rk(I)=\rk(B)+\rk(I_\sigma)$, and we claim
$\rk(I_\sigma)=\rk(B)+\rk(C(\sigma))$.

For any element $i\in I_\sigma$ we have $i=yj$ for some $y\in B\sigma$ and
$j\in I$, and hence $i$ inverts some $y\in B\sigma$.
This element $y$ is {\em unique:} if $i$ inverts $y$ and $by$ with $b\in
B^\#$, then $y^{-1}b^{-1}=(by)^i=b^iy^{-1}$ and thus $b^{-1}\in
B^{i y^{-1}}$, $i y^{-1}\in N(B)$, and $i\in N(B)$, contradicting the
choice of $i\in I_1$. So we have a definable function
$$\beta:I_\sigma\to B\sigma$$
defined by $\beta(i)^i=(\beta(i))^{-1}$. 

As $B\sigma$ is a single conjugacy class under the
action of $B$, the rank of the fibers $\beta^{-1}(y)$ is a constant $f$, 
and hence
$\rk(I_\sigma)=\rk(B\sigma)+f=\rk(B)+\rk(\beta^{-1}(\sigma))$. 

It suffices therefore to show
$\rk(\beta^{-1}(\sigma))=\rk(C(\sigma))$.  
Fix $i$ in $\beta^{-1}(\sigma)$.
Then $i=\sigma i'$ with $i'\in I$ inverting $\sigma$. 

We claim first
$$iC\oo(\sigma)\includedin \beta^{-1}(\sigma)\leqno(*)$$
Observe first that $C(\sigma)$ contains no involutions, by
Fact \ref{BN:stronglyreal}, bearing in mind Lemma 4.12.

For $g\in C\oo(\sigma)$ and $j=ig$, we have $j\in I$ by Lemma
\ref{inversion}. 
Furthermore $j=\sigma\cdot i'g$, and $i'g\in I$ by Lemma \ref{inversion}.
So $j\in B\sigma I\intersect I$.
By construction $j$ inverts $\sigma$.  It remains to check that
$j\notin M$.
If $j\in M$ then $j\in B$, and since $j$
inverts $\sigma$, and $\sigma$ is of odd order modulo $B$, 
this gives a contradiction. So $j\in I_1$ 
and $\beta(j)=\sigma$.
This gives $(*)$, and in particular
$\rk(\beta^{-1}(\sigma)) \ge \rk(C\oo(\sigma))=\rk(C(\sigma))$.

Conversely, for $j\in \beta^{-1}(\sigma)$, since $j$ inverts $\sigma$
we have $ij\in C(\sigma)$, that  is $i\cdot
\beta^{-1}(\sigma)\includedin C(\sigma)$, and thus
$\rk(\beta^{-1}(\sigma))\le \rk(C(\sigma))$. Our claim follows.
\end{proof}

With this our analysis is complete. $BC(\sigma)B$ and $BI_1$ are
disjoint generic subsets of $G$ by Lemmas \ref{1stgeneric}, \ref{rankBCB}
and \ref{rankG}, and \ref{disjoint}, which is a contradiction.
Thus Case I cannot occur in Pr\"ufer $2$-rank two or more.


%% file: min04.tex

\section{Case II}
\label{caseII}

The remaining case will require a longer analysis, and some more
theoretical preparation. We begin afresh.

Our standing hypotheses and notation are as follows.

\begin{notation}
\ 
\begin{enumerate}
\item $G$ is a minimal connected simple group of finite Morley rank and of
odd type.
\item $B$ is a Borel subgroup of $G$.
\item $M=N(B)$ is strongly embedded in $G$.
\item[$(\dag)$]$G$ has Pr\"ufer $2$-rank at least two.
\end{enumerate}
\end{notation}

The operative assumption for the remainder of the article
will be as follows.
$$\hbox{For $i$ an involution of $B$, we have $C\oo(i)<B$.}
\leqno(\hbox{\bf Case II})$$

It will be convenient to state this in a slightly stronger form.

\begin{lemma}\label{F(B)}
$F(B)$ contains no involutions.
\end{lemma}
\begin{proof}
Let $S$ be a Sylow 2-subgroup of $F(B)$. Then $A=\Omega_1(Z(S))$ 
is characteristic in $F(B)$ and normal in $B$. 
As $G$ is of odd type, the group $A$ 
is finite as well as $B$-invariant, and hence central in $B$.
By our case assumption (II), $A=1$ and hence $S=1$.
\end{proof}

\subsection{Tameness and Fact \protect\ref{jaligotlemma}}\label{tamecase}

Let us first indicate why Case II disappears quickly 
if we assume tameness.
First, by an easy argument (Lemma \ref{I*} below) our case assumption
produces an involution $w$ outside $M$ such that $B\intersect B^w$
is infinite. By strong embedding $M\intersect M^w$ contains no
involutions and hence the same applies to $H=(B\intersect B^w)\oo$. 
It then follows by an application of tameness that 
$H$ is contained in the Fitting subgroup of both $B$ and $B^w$, and
this contradicts Fact \ref{jaligotlemma}.

We must make a distinction between the two applications of tameness in the
foregoing argument. As mentioned earlier, we will make use of a weak
analog of Fact \ref{jaligotlemma} which does not require
tameness. However the claim that a connected group without involutions
must lie in the Fitting subgroup is a direct application of tameness
with no obvious analog in general. So before entering into the general
case, let us first give an argument which uses only Fact
\ref{jaligotlemma} and therefore serves as a reasonable point of
departure for the general case. This will necessarily be more
elaborate than the argument just given, and will only be sketched,
since in any case it will need to be redone afterward at a greater
level of generality.

We begin again with an involution $w\notin M$ for which $B\intersect
B^w$ is infinite, and more particularly 
$$T(w)=\{a\in B:a^w=a^{-1}\}$$
is infinite. Observe that $T(w)\includedin B\intersect B^w$, and under
the assumption that Fact \ref{jaligotlemma} applies, the group
$B\intersect B^w$ is abelian (as the natural map $B\intersect B^w\to
B/F(B)\times B^w/F(B^w)$ is injective). So $T(w)$ is a group in this
case. Let $H=(B\intersect B^w)\oo$ and consider a maximal definable 
connected subgroup $\hat H$ containing $H$ of the form $(B\intersect
B_1)\oo$, 
with $B_1\ne B$ a Borel, not necessarily standard. Again, Fact
\ref{jaligotlemma} implies that $\hat H$ is abelian. Note that $\hat H\le
N\oo(H)$ and $N\oo(H)$ is $w$-invariant. We will take $B_1\ge
N\oo(H)$, and with some effort we may make $B_1$ $w$-invariant as
well. 

Observe that $\hat H$ cannot be a Carter subgroup of $B$: otherwise, it
contains a Sylow 2-subgroup of $B$ by Lemma \ref{carter-sylow},
forcing $B_1=B$. By the maximality of $\hat H$ it then follows that
$\hat H$ is a Carter subgroup of $B_1$. On the other hand by a Frattini
argument $w$ normalizes a Carter subgroup of $B_1$ and hence some
conjugate $w_1$ of $w$ (under the action of $B_1$) normalizes $\hat H$.
As $w_1$ normalizes $N\oo(\hat H)$, it follows easily that $w_1$ normalizes,
and hence lies in, the group $B$. This is the fundamental setup that
will be reached below, in general.
$$\hbox{$w_1\in I(B)$ normalizes $B_1$, and is conjugate to $w$ under
  the action of $B_1$.}$$

We now use the characteristic zero unipotence theory 
introduced in \cite{Bu-Th,Bu-SF}, 
and in particular the notation $U_0(T(w))$, 
which represents a kind of unipotent radical,
and the notion of ``reduced rank''. We also use 
the associated generalized Sylow theory of \cite{Bu-Th,Bu-GCS}.
(This machinery will be presented in more detail below.)

The group $T(w)$ is a nontrivial definable abelian group, inverted by
$w$, and it can be shown to be torsion free, hence connected.  Since
$T(w)$ is nontrivial, solvable, and torsion free, it has a nontrivial
``0-unipotent radical'' $T_w:=U_0(T(w))$. Let the maximal reduced rank
associated with $T_w$ be $r$; then $T_w=U_{0,r}(T(w))$.  One can
extend $T_w$ to a $w$-invariant Sylow $U_{0,r}$-subgroup $P_r$ of
$B_1$, and any two such are conjugate in $C_{B_1}(w)$ (Lemma
\ref{invariantsylow}, below). As $P_r$ 
contains a nontrivial $U_{0,r}$-subgroup inverted by $w$, the same
applies to every Sylow $U_{0,r}$-subgroup of $B_1$ which is normalized
by $w$.  Since $w$ and $w_1$ are conjugate under the action of $B_1$,
the involution $w_1$ satisfies an analogous condition: any
$w_1$-invariant Sylow $U_{0,r}$-subgroup of $B_1$ contains a
nontrivial $U_{0,r}$-subgroup inverted by $w_1$.

On the other hand, if $Q_r=U_{0,r}(\hat H)$, then $N_B\oo(Q_r)\ge
N_B\oo(\hat H)>\hat H$ and hence by maximality $N\oo(Q_r)\le B$,
from which it follows that $Q_r$ is a Sylow $U_{0,r}$-subgroup of
$B_1$, and of course $w_1$-invariant. So there must be a nontrivial
$U_{0,r}$-subgroup $A_r$ of $\hat H$ inverted by $w_1$. 

Now one can show easily that for any $s$, $F_s(B_1):= U_{0,s}(F(B_1))$ 
is either
contained in $\hat H$, or meets $\hat H$ trivially (cf.~the proof of Lemma 
\ref{U0rhatH-}, claim (*), below). In either case,
$A_r$ commutes with $F_s(B_1)$: if $F_s(B_1)\le \hat H$, 
this holds because
$\hat H$ is abelian, and if $F_s(B_1)\intersect \hat H=1$, then
$w_1$ inverts $F_s(B_1)$, and consideration of the action of $w_1$
on the group $F_s(B_1)A_r$ leads to the desired conclusion.
From all of this it follows that
$A_r$ centralizes $F(B_1)$, and hence lies
in $F(B_1)$. But $A_r=[w_1,A_r]\le F(B)$ as well, and by Fact
\ref{jaligotlemma} we find $A_r=1$, a contradiction.

We will argue in the remainder of the paper that some of the
applications of Fact \ref{jaligotlemma} can be avoided and others
replaced by a more complicated, but general, form of that result. One
major alternative that arises in general is that the group $\hat H$ in
question is nonabelian, but this produces a rather well defined
configuration which turns out to have a good deal in common with the
abelian case.  In fact as we will see the analysis follows much the
same lines whether $\hat H$ is abelian or nonabelian.

\subsection{The setup}

Now dropping the assumption that Fact \ref{jaligotlemma} applies and
returning to the general case, we begin by showing that the setup with
which we started our analysis above can in fact be reached.

\begin{notation}
\ 
\begin{enumerate}
\item $I=I(G)$.
\item For $w\in I$, let $T[w]$ be $\{a\in B:a^w=a^{-1}\}$.
\item Let $I^*=\{w\in I\setminus N(B):\rk(T[w])\ge \rk (I(B))\}$
\end{enumerate}
\end{notation}

The ungainly notation $T[w]$ is intended to reflect the
ungainly nature of the set involved, which in general need
not be a group.
Our main  concern is that $I^*$ should be nonempty, and this is
afforded by Lemma \ref{I*} below.

We use the following general fact.

\begin{fact}[{\cite[2.36]{CJ01}}]\label{conjugacy}
Let $G$ be a connected simple group of finite Morley rank, let $M$ be a
proper definable subgroup of $G$, and let $X$ be a conjugacy class
in $G$. Then $\rk(X\intersect M)<\rk(X)$.
\end{fact}

\begin{lemma}\label{I*}
$I^*$ is generic in $I$.
\end{lemma}
\begin{proof}
Let $i\in I(B)$. Then we have $\rk(I)=\rk(G)-\rk(C(i))$ and thus
$\rk(I(B))={\rk(B)- \rk(C(i))}=\rk(I)-\rk(G/B)$.

By Fact \ref{conjugacy}, $I\setminus N(B)$ is generic in $I$.
It will suffice to prove that the set $I'$ defined as
$$\{w\in I\setminus B:\rk(T[w])<\rk(I(B))\}$$
is nongeneric in $I$.

Now $T[w]=w\cdot(wB\intersect I)$, so $\rk(T[w])=\rk (wB\intersect I)$. 
For $w\in I'$, it follows that $\rk(wB\intersect I)<\rk(I(B))$.
Hence $\rk(I'\intersect X)<\rk(I(B))$ for $X$ any left coset of $B$ in $G$.
As $I'=\Union_{X\in G/B}(I'\intersect X)$, we find
$$\rk(I')<\rk(G/B)+\rk(I(B))=\rk(I),$$
as claimed.
\end{proof}

Now we can fix our notation for the remainder of the argument.

\begin{notation}\label{caseIIhypotheses}
\ 
\begin{enumerate}
\item Fix $w\in I^*$. Set $H=(B\intersect B^w)\oo$.
\item Let $B_1\ne B$ be a Borel subgroup containing $H$, chosen so as to
maximize $(B\intersect B_1)\oo$.
\item Let $\hat H=(B\intersect B_1)\oo$.
\end{enumerate}
\end{notation}

There is some latitude in the choice of $\hat H$, and with $\hat H$ 
fixed there
is some latitude in the choice of $B_1$, which will be examined
more closely subsequently. Our initial goal is to show that $B_1$ can
be chosen to be $w$-invariant. Along the way we will acquire other
useful information. 


%% file: min05.tex

\section{$w$-Invariance}\label{sec:invariance}

We begin our analysis of Case II. We have the notations $G,B,w, H,
\hat H,B_1$ as laid out in the previous section (and also $T[w]$,
which will be needed at a later stage, when the configuration is
clearer). 
In particular $H\le \hat H=(B\intersect B_1)\oo$ with $B_1$ a Borel
subgroup. As always, the Pr\"ufer 2-rank is assumed to be at least two.

Our goal in the present section is to show that $B_1$ can be chosen
to be $w$-invariant. 

\subsection{Unipotence Theory}

We use the $0$-unipotence theory of \cite{Bu-Th}
(cf.~\cite{Bu-SF,FJ-CS}). We use the notation $\rred(K)$ for the {\em
maximal reduced rank} of a group of finite Morley rank $K$,
which is the largest integer $r$ for which $U_{0,r}(K)\ne 1$ (or $0$
if there is no such $r$). 
If $K$ is a group of finite Morley
rank then we set $U_0(K)=U_{0,\rred(K)}(K)$. 

We make use of the following from the general $0$-unipotence theory.


\begin{fact}\label{unipotencebasics}
Let $K$ be a connected solvable group of finite Morley rank.
\begin{enumerate}
\item $U_0(K)\le F(K)$
\item \label{ub2}
If $K$ is nilpotent and $K=U_{0,r}(K)$, then $K'=U_{0,r}(K')$.
\item \label{ub5}
If $1\to K\to H\to \bar H\to 1$ is a short exact sequence of definable
groups with $U_{0,r}(K)=K$ and $U_{0,r}(\bar H)=\bar H$,
then $U_{0,r}(H)=H$. Conversely, if $U_{0,r}(H)=H$ then $U_{0,r}(\bar
H)=\bar H$.
\item \label{ub3}
$\rred(K)=\rred(Z(F(K)))$.
\item \label{ub4}
If $K=U_0(K)$ and $r=\rred(K)$, then for any proper definable
  subgroup $K_0< K$ we have $U_{0,r}(N_K(K_0))>K_0$.
\item \label{nilstructure}
  If $K$ is nilpotent, then $K=B*T$ with $B,T$ definable, $B$ of
  bounded exponent, and $T$ divisible; $B$ is the 
  central product of the finitely many subgroups $U_p(K)$,
  with $p$ prime, for which
   $U_p(K)\ne 1$; $T$ is the central product of the finitely many 
  groups $U_{0,s}(K)$
  for which $U_{0,s}(K)\ne 1$, together with the group
  $K_\infty=d(T\tor)$, the definable closure of the torsion subgroup
  of $T$.
\end{enumerate}
\end{fact}

The first three points are given in \cite[2.16, 2.17, 2.11]{Bu-SF},
and the fourth follows from the first two. The last two are given
in \cite[Lemma 2.4, Corollary 3.5]{Bu-GCS}. 
These points are also in \cite{Bu-Th} as
Theorem 2.21, Lemmas 2.23, 2.12, 2.26, and 2.28, and Theorem 2.31,
respectively. 

We note that the first five points are close analogs of more elementary 
properties of $U_p$ for $p$ prime, which we will use without special
comment. 

The next result has a similar character, but is used less often.

\begin{fact}[{\cite[4.9]{Bu-Th},\cite[Lemma 4.4]{Bu-GCS}}]
\label{UbyU}
Let $G=H_1\cdot H_2$ be a group of finite Morley rank with
$H_i$ a definable nilpotent $U_{0,r_i}$-subgroup of $G$ for $i=1$ or $2$, 
where $H_1$ is normal in $G$, and $r_2\ge r_1$.
Then $G$ is nilpotent.
\end{fact}

In the case that interests us, $r_1=r_2$.

The utility of this abstract theory of unipotence
often depends on the following.

\begin{fact}[{\cite[2.15]{Bu-SF},\cite[2.19]{Bu-Th}}]
\label{completeness} 
Let $H$ be a connected solvable group of finite Morley rank such that
$U_p(H)=1$ for all primes, as well as for $p=0$. Then $H$ is a good
torus; that is, a divisible abelian group in which every definable
subgroup is the definable closure of its torsion subgroup.
\end{fact}

The following generation principle can be very useful.

\begin{fact}[{\cite[2.41]{Bu-Th},\cite[Theorem 2.9]{Bu-GCS}}]
\label{U0r-generation} 
If a nilpotent $U_{0,r}$-group is generated by a family $\FF$
of definable subgroups, then it is generated by the 
the family $\FF_r=\{U_{0,r}(X):X\in \FF\}$.
\end{fact}

We have alluded also to a Sylow theory. By definition a {\em Sylow
$U_{0,r}$-subgroup} is a maximal definable nilpotent
$U_{0,r}$-subgroup. The conjugacy theorem applies, at least in a
solvable context.

\begin{fact}[{\cite[4.16, 4.18]{Bu-Th},\cite[Theorem 6.5]{Bu-GCS}}]
\label{U0rconjugacy}
Let $K$ be a solvable group of finite Morley rank. Then for each $r$,
its Sylow $U_{0,r}$-subgroups are conjugate.
\end{fact}

Some further connections between the Sylow theory and the Carter
theory will be recalled when needed.

\subsection{A Uniqueness Lemma}

As we have noticed previously, when Fact \ref{jaligotlemma} applies,
the intersections of distinct Borel subgroups are abelian. 
The following is an indication of the tension that arises when the
latter condition fails. We give it in a general setting.

\begin{lemma}[{\cite[Theorem 4.3]{Bu-JL}}]\label{uniqueness-1}
Let $G$ be a minimal connected simple group of finite Morley rank,
let $B_1$, $B_2$ be distinct Borel subgroups of $G$, and
$H=(B_1\intersect B_2)\oo$. Suppose that $H$ is nonabelian. Then the
following conditions are equivalent.
\begin{enumerate}
\item $B_1$ and $B_2$ are the only Borel subgroups of $G$ containing
  $H$.
\item If $B_3$ and $B_4$ are distinct Borel subgroups containing $H$,
then $(B_3\intersect B_4)\oo=H$.
\item If $B_3\ne B_1$ is a Borel subgroup containing $H$, then
  $(B_1\intersect B_3)\oo=H$. 
\item $C\oo(H')$ is contained in $B_1$ or $B_2$.
\item $\rred(B_1)\ne \rred(B_2)$
\end{enumerate}
\end{lemma}

Note that clauses $(2,3)$ express the maximality of $H$ in two
different senses: in $(2)$ we vary $B_1$ and $B_2$, while in $(3)$ we
hold $B_1$  fixed and only vary $B_2$. The first clause is an even
more extreme form of maximality, while the last two clauses 
provide remarkably simple criteria for identifying such pairs
$B_1,B_2$. They both follow fairly readily from clause (2), but the
converse is more subtle. 

We note the following consequence.

\begin{lemma}\label{uniqueness-2}
Let $G$ be a minimal connected simple group of finite Morley rank,
and $H$ a proper connected definable nonabelian
subgroup of $G$. 
Then $HC\oo(H')$ is contained in a unique Borel subgroup
of $G$, and in particular $N\oo(H')$ is contained in a unique Borel
subgroup of $G$.
\end{lemma}
\begin{proof}
If $B_1,B_2$ are distinct Borel subgroups of $G$ containing
$HC\oo(H')$, let $H_1=(B_1\intersect B_2\oo)$. 
Then $H_1C\oo(H_1')$ is
also contained in $B_1$ and $B_2$. So by the previous result, 
$\rred(B_1)\ne \rred(B_2)$; we may suppose $\rred(B_1)>\rred(B_2)$.
Then as $H_1'\le F(B_1)$ \cite[Cor.~9.9]{BN}, we have
$U_0(B_1)\le C(H_1')$ by 
Fact \ref{unipotencebasics} (\ref{nilstructure}), and hence
$U_0(B_1)\le B_2$, contradicting $\rred(B_1)>\rred(B_2)$.
\end{proof}

Now let us return to the case at hand, in which the pair
$(B,B_1)$ plays the role of $(B_1,B_2)$ above, so that our $\hat H$ 
corresponds to the $H$ of our lemma. Then the third clause expresses
exactly the maximality condition 
that we have imposed on our pair $(B,B_1)$, and
hence when $\hat H$ is nonabelian all of these conditions apply.
This fact will play a leading role throughout the rest of the
analysis, coming into play whenever $\hat H$ is assumed nonabelian.
The main point here is that we find ourselves in the situation
described by the third clause, while the second clause is the one
most conveniently adopted as a point of departure for the detailed
analysis undertaken in \cite{Bu-JL}; since we need this only when
$\hat H$ is nonabelian, we can cite \cite{Bu-JL} freely in such cases.
For the record, we state this in the slightly stronger form afforded
by Lemma 6.6, part (1).

\begin{corollary}
If $\hat H$ is nonabelian, then $B$ and $B_1$ are the only Borel
subgroups containing $\hat H$.
\end{corollary}

Furthermore, applying Lemmas \ref{uniqueness-1} and \ref{uniqueness-2},
we now find the desired $w$-invariant Borel subgroup
in one important case.

\begin{lemma}\label{Hnonabelian}
If $H$ is nonabelian, then $B_1$ can be chosen to be $w$-invariant.
\end{lemma}
\begin{proof}
Take $\tilde B_1$ to be a Borel subgroup containing
$N\oo(H')$. By Lemma \ref{uniqueness-2}, $\tilde B_1$ is the only such
Borel subgroup, and is therefore $w$-invariant, and in particular
$\tilde B_1$ is not $B$.
Let $\tilde H= (B\intersect \tilde B_1)\oo$.  Then 
$C\oo({\tilde H}')\le C\oo(H')\le \tilde B_1$.
So by Lemma
\ref{uniqueness-1}, the groups $B$ and $\tilde B_1$ are the only Borel
subgroups containing $\tilde H$, and in particular $\tilde H$ is
maximal.  All of our conditions are met with $B_1=\tilde B_1$.
\end{proof}

\subsection{Extension Lemma}

We recall that the notation and operative hypotheses were
established at the end of the previous section in Notation
\ref{caseIIhypotheses}. We will insist somewhat on this point for the
remainder of the present section, because some of the work takes place
at a sufficient level of generality to allow for its reuse in the next
section, and goes beyond the immediate needs of the moment. We will of
course have to track carefully what additional hypotheses are imposed
in particular results.

\begin{lemma}
With our current hypotheses and notations concerning $G$, $B$, and $w$,
suppose $K\le G$ is a maximal proper
definable connected $w$-invariant subgroup of $G$.  
Then $K$ is nonabelian.
\end{lemma}
\begin{proof}
Suppose on the contrary that $K$ is abelian, and let $B_w$ be the
conjugate of $B$ containing $w$.

By maximality of $K$, it follows that 
$N\oo(K_0)=K$ for any nontrivial $w$-invariant
definable subgroup $K_0$ of $K$.  

In particular, if 
$K_0=(K\intersect B_w)\oo\ne 1$, 
then $N\oo(K_0)=K$.
Then $N_{B_w}\oo(K_0)=K_0$, so $K_0$
is a Carter subgroup of $B_w$, in view of Fact \ref{nearcarter}.  
Then $K_0$ contains a Sylow $2$-subgroup
$S_w$ of $B_w$ by Lemma \ref{carter-sylow}.  It follows from Lemma 
\ref{blackhole} that $K\le B_w$, and hence 
by maximality that $K=B_w$. Thus $B_w$ is abelian, so $B$ is abelian,
contradicting our current case hypothesis, namely $C\oo(i)<B$.

So suppose that $K\intersect B_w$ is finite, and in particular
$C_K(w)$ is finite. Then $w$ inverts $K$ \cite[p.~78, Ex.~13]{BN}. 
In particular,
every subgroup of $K$ is $w$-invariant. It follows that distinct
conjugates of $K$ are disjoint: if $K\intersect K^g\ne 1$, then
$N\oo(K\intersect K^g)$ is equal to both $K$ and $K^g$, 
and $K=K^g$.

As $K$ has finite index in its normalizer, and has pairwise disjoint
conjugates, the union $\Union_{g\in G}K^g$ is generic in $G$ by Fact
\ref{genericcovering}, and hence meets $\Union_{g\in G}B^g$ 
nontrivially by Lemma
\ref{Bgenericcovering}. Let $B_0$
be a standard Borel subgroup that meets $K$ nontrivially.

Now $K\intersect B_0$ is $w$-invariant and nontrivial, and hence
$N\oo(K\intersect B_0)\le K$. It follows by Fact \ref{nearcarter}
that $K\intersect B_0$ is a Carter subgroup of $B_0$, 
and hence contains a Sylow 2-subgroup of $B_0$. 
It then follows from Lemma  \ref{blackhole} 
that $K\le B_0$. Then again, $B_0=K$ should be abelian, contradicting 
our case assumption.
\end{proof}

\begin{corollary}\label{invariantBorel}
Under the assumptions of the present section,
if $K$ is a nontrivial proper 
connected definable $w$-invariant subgroup of
$G$, then $K$ is contained in a $w$-invariant Borel subgroup of $G$.
\end{corollary}
\begin{proof}
We may take $K$ to be maximal under the stated conditions, and hence
nonabelian by the preceding lemma. Then by Lemma \ref{uniqueness-2},
$N\oo(K')$ is contained in a unique Borel subgroup of $G$, which is
again $w$-invariant.
\end{proof}

Using this result, we deal easily with the case in which $\hat H$ is
abelian. 

\begin{lemma}\label{Hhat-abelian}
Suppose that $\hat H$ is abelian. Then $B_1$ can be chosen to be
$w$-invariant.
\end{lemma}
\begin{proof}
We have $\hat H\le N\oo(H)$ and the latter is $w$-invariant,
So $N\oo(H)$ can be extended to a $w$-invariant Borel subgroup 
$B_1$.  So $B_1\ge \hat H$, 
and as $B_1$ is $w$-invariant, we have $B_1\ne B$.
\end{proof}

\subsection{Maximal pairs}

We have seen that $B_1$ can be chosen $w$-invariant if $\hat H$ is
abelian, or if $H$ is nonabelian, and we now consider the remaining
possibility: $H$ is abelian, while $\hat H$ is nonabelian. As we have
seen in Lemma \ref{uniqueness-1}, 
if $\hat H$ is fixed then the Borel subgroup $B_1$ is uniquely
determined by the conditions $\hat H\le B_1\ne B$. However, there is
still some latitude in the choice of $\hat H$.

According to Fact \ref{jaligotlemma}, this situation cannot arise in
the tame case. More generally,  while not visibly
contradictory, this configuration is tightly constrained in general,
as described in \cite[9.2]{Bu-Th} and \cite{Bu-JL}.

\begin{definition}
Let $G$ be a group of finite Morley rank, and $B,B_1$ two Borel
subgroups of $G$. We call the pair $(B,B_1)$ a {\em maximal pair}
if $(B\intersect B_1)\oo$ is maximal, among all connected components
of intersections of distinct Borel subgroups of $G$;
in other words, clause (2) of Lemma \ref{uniqueness-1} applies.
\end{definition}

When the intersection in question is nonabelian, we need not be very
particular about the notion of maximality invoked, but the definition
adopted here coincides with the one given in \cite{Bu-JL}, which is well
suited also to analysis in some abelian cases.

In the first place we have the following. We will only apply this when
our group $\hat H$ is nonabelian, but we take note of the slightly
greater generality achieved in \cite{Bu-JL}.

\begin{fact}[{\cite[9.2 (2, 4d)]{Bu-Th},\cite[Theorem 4.5, (1,5)]{Bu-JL}}]
\label{9.2Carter}
Let $G$ be a minimal connected simple group of finite Morley rank,
and $(B,B_1)$ a maximal pair of Borel subgroups of $G$. Let
$H=(B\intersect B_1)\oo$ and let $Q$ be a Carter subgroup of $H$.
Suppose that $H$ is nonabelian, or more generally 
that $F(B_1)\intersect F(B_2)$ is nontrivial.
Then 
\begin{enumerate}
\item $\rred(B)\ne \rred(B_1)$;
\item If $\rred(B)<\rred(B_1)$, then 
\begin{enumerate}
\item$\rred(H)=\rred(B)$ and
\item  $Q$ is a Carter subgroup of $B_1$. 
\end{enumerate}
\end{enumerate}
\end{fact}

Let us apply this now to the case at hand, in which our pair $(B,B_1)$
plays the role of the given pair $(B,B_1)$, 
and we suppose that $\hat H$, which plays the role of $H$, 
is nonabelian, so that we do indeed have a maximal pair
in view of Lemma \ref{uniqueness-1}. 
Let us make the relevant conclusions
explicit in this case. Note that one hypothesis of the foregoing 
fact becomes a conclusion in our context.

\begin{lemma}\label{9.2Carter'}
Under our present hypotheses and with our present notation, if 
$Q$ is a Carter subgroup of $\hat H$, then $Q$ is not a Carter
subgroup of $B$. If in addition
$\hat H$ is nonabelian, 
then
\begin{enumerate}
\item $\rred(B)<\rred(B_1)$.
\item $Q$ is a Carter subgroup of $B_1$.
\end{enumerate}
\end{lemma}
\begin{proof}
If $Q$ is a Carter subgroup of $B$, then by Lemma \ref{carter-sylow},
the group $Q$ contains a Sylow $2$-subgroup $S$ of $B$. Thus $S\le B_1$,
and by Lemma \ref{blackhole} we find $B_1=B$, a contradiction. This
proves the first point.

Now take $\hat H$ to be nonabelian.
Then by Fact \ref{9.2Carter}, the group $Q$ is a Carter subgroup of
whichever group, $B$ or $B_1$, has the larger reduced rank.
As this cannot be $B$, it must be $B_1$, and our claims follow.
\end{proof}

The Carter subgroups of $\hat H$ play a central role in what follows,
largely because of Lemma \ref{QandB} below, which treats both the case
in which $\hat H$ is abelian and in which it is not. For the
nonabelian case, we will need some additional information from
\cite{Bu-JL}, particularly bearing on the commutator subgroup $\hat
H'$ (Fact \ref{JL:H'} below).

\begin{definition}
Let $H$ be a solvable group of finite Morley rank. Then $H$ is {\em
  rank-homogeneous} if it satisfies the following conditions.
\begin{enumerate}
\item $H$ is torsion free.
\item For $r<\rred(H)$, $U_{0,r}(H)=1$.
\end{enumerate}
\end{definition}

Note that these two clauses force $H=U_0(H)$, and hence $H$ is
nilpotent by Fact \ref{unipotencebasics} (1), which
clarifies the meaning of the definition, particularly if the
structure theory of Fact \ref{unipotencebasics} (6) is kept in mind. 

\begin{notation}
Let $H$ be a group of finite Morley rank, and $r\ge 0$. Then
$F_r(H)$ denotes $U_{0,r}(F(H))$. (One prefers $r>0$: $F_0(H)=1$.)
\end{notation}

\begin{fact}[{\cite[9.2 (5a,1,5b,3,5d,5c)]{Bu-Th},
\cite{Bu-JL}}]
\label{JL:H'}
Let $G$ be a minimal connected simple group of finite Morley rank,
and $(B,B_1)$ a maximal pair of Borel subgroups of $G$. Let
$H=(B\intersect B_1)\oo$ and let $Q$ be a Carter subgroup of $G$.
Suppose that $H$ is nonabelian, and let the notation be chosen so that
$\rred(B)<\rred(B_1)$. 
Then we have the following.
\begin{enumerate}
\item $H'$ is rank-homogeneous.
\item If $r'=\rred(H')$, then $U_{0,r'}(H)=F_{r'}(H)$.
\item $F(H)$ is abelian.
\item $Q$ is abelian.
\item $Q_{r'}\le Z(H)$ is nontrivial.
\item For any nontrivial definable  subgroup $X$ of $H$ which
  is contained in $H'$,  we have $N\oo(X)\le B_1$.
\item \label{H':5}$F_{r'}(H)<F_{r'}(B)$.
\end{enumerate}
\end{fact}

The first four points
are covered in \cite{Bu-JL} with the encyclopedic Theorem 4.5.
More precisely: since $H'\le F(B_1)\intersect F(B_2)$, 
$(1)$ is contained in 4.5 $(6)$; $(2)$ is contained in 4.5 $(3)$;
both $(3)$ and $(4)$ are contained in 4.5 $(2)$; 

The fifth and sixth points are given as Lemma 3.23 and Corollary 3.29
of \cite{Bu-JL} respectively. The last point combines part of 
4.5 $(2)$ with Lemma 3.13 of \cite{Bu-JL}.

This fact also generalizes to the case in which the group
$F(B_1)\intersect F(B)$, which may well be larger than $H'$,
is nontrivial, and one can replace $H'$ by that larger group for these
purposes. We do not need this refinement.

We proceed now with our analysis. While we would be free at this point
to assume $\hat H$ is nonabelian, we need some of these results in the
following section in a broader context, so any assumptions needed on
$H$ or $\hat H$ will be stated explicitly as required; and
otherwise we will take note of their absence.

\begin{lemma}\label{QandB}
Under our present hypotheses and notations, but
without additional assumptions on $H$ or $\hat H$,
let $Q$ be a Carter subgroup of $\hat H$, and for any $r$ let
$Q_r=U_{0,r}(Q)$.
 Then the following hold.
\begin{enumerate}
\item $N(Q)\le N(B)$. 
\item If $\hat H$ is nilpotent and $Q_r>1$, then $N(Q_r)\le N(B)$.
\end{enumerate}
\end{lemma}
\begin{proof}
First we claim
$$N_B\oo(Q)\not \le \hat H\leqno(*)$$
If $N_B\oo(Q)\le \hat H$ then $N_B\oo(Q)\le N_{\hat H}(Q)=Q$ and thus
$Q$ is a Carter subgroup of $B$ by Fact \ref{nearcarter},
contradicting Lemma \ref{9.2Carter'}. So $(*)$ holds.

Now we divide into two cases, according as $\hat H$ is or is not
nilpotent. 

Suppose first
$$\hbox{$\hat H$ is nilpotent}\leqno(1)$$
Then $Q=\hat H$ and by $(*)$ we have $N_B\oo(Q)>Q$. By maximality 
of $\hat H$, the group $B$ is the only Borel subgroup containing
$N_B\oo(Q)$. Hence $B$ contains $N\oo(Q)$, and is the only Borel
subgroup containing $N\oo(Q)$. 
From the last point it follows that $N(Q)$ normalizes $B$.
Now if $Q_r>1$, then
as $N\oo(Q)\le N(Q_r)<G$, it also follows that $B$ is the unique Borel
subgroup containing $N\oo(Q_r)$. Hence $N(Q_r)\le N(B)$.

Now suppose on the contrary
$$\hbox{$\hat H$ is nonnilpotent}\leqno(2)$$
and in particular nonabelian.
We set $\hat r'=\rred(\hat H')$.

Let $Q_{\hat r'}=U_{0,\hat r'}(Q)$.  By Fact \ref{JL:H'} (5) we have
$Q_{\hat r'}\le Z(\hat H)$, so $\hat H\le N(Q_{\hat r'})$, and since
$N\oo(Q)\le N\oo(Q_{\hat r'})$ we find $\hat H<N_B\oo(Q_{\hat r'})$.
Since $Q_{\hat r'}$ is nontrivial by Fact \ref{JL:H'} (5), 
it follows by maximality of $\hat H$ 
that the group $N\oo(Q_{\hat r'})$ is contained in $B$ and in no other
Borel subgroup.  In particular as $N(Q)$ normalizes $Q_{\hat r'}$,
it follows that $N(Q)$ normalizes $B$.
\end{proof}

Recall now that $H=(B\intersect B^w)\oo$. We insert a lemma
which simplifies matters somewhat as far as $H$ is concerned. 

\begin{lemma}\label{torsionfree}
Under our present hypotheses and notations, but without requiring
$\hat H$ to be nonabelian, the intersection
$B\intersect B^w$ is torsion free.
\end{lemma}
\begin{proof}
Suppose on the contrary $x\in B\intersect B^w$ has prime order $p$.
Let $P\le B$ be a Sylow\oo\ $p$-subgroup of $B$ containing $x$.  As
$B$ is solvable, the group $P$ is locally finite and we can use the
structure theory of \cite[6.20]{BN}; since in addition
\cite[9.39]{BN} $P$ is connected, it follows that either
$U_p(P)>1$, or $P$ is a $P$-torus.

If $U_p(B)\ne 1$, then $U_p(C_{U_p(B)}\oo(x))$ is a nontrivial
$p$-unipotent group and hence is contained in a unique Borel subgroup
of $G$ (Lemma \ref{punipotent}), which must be $B$. So $C\oo(x)\le
B$. But if $U_p(B)\ne 1$ then $U_p(B^w)\ne 1$, so similarly
$C\oo(x)\le B^w$, and now the uniqueness statement yields $B=B^w$, a
contradiction.

So $U_p(B)=1$ and $P$ is a $p$-torus. Then by Lemma
\ref{carter-sylow}, $P$ is contained in a Carter subgroup $R$ of $B$,
and by the same lemma 
$R$ also contains a Sylow 2-subgroup $S$ of $B$. So as $R$ is
nilpotent, $P$ and $S$
commute, and it follows that $C(x)$ contains $S$. Similarly $C(x)$
contains a Sylow 2-subgroup $S_1$ of $B^w$. Then $\Omega_1(S)$
normalizes $C\oo(x)$ and hence by Lemma \ref{blackhole}, we have
$C\oo(x)\le B$ and $S_1\le B$, forcing $B=B^w$, a contradiction.
\end{proof}

In consequence we have $H=B\intersect B^w$, $H$ is torsion free, and
$T[w]\includedin  H$.
When $H$ is abelian, it follows that $T[w]$ is a subgroup of $H$
inverted by $w$. Note also that $T[w]$ contains some infinite
definable abelian subgroups inverted by $w$, which is a small 
start in the right direction. Eventually we will arrive at the case in
which $T[w]$ is itself an abelian group.

\subsection{Invariance of $B_1$}

Throughout this subsection we deal with the case in which $\hat H$ is
nonabelian and $H$ is abelian, though our preparatory work is more general. 
Recall that we still have some latitude in the choice of $\hat H$.

\begin{notation}
\ 
\begin{enumerate}
\item $\hat r'=\rred(\hat H')$.
\item Set $H^-=d(T[w])$, the smallest definable subgroup containing $T[w]$.
\item Set $r^-=\rred(H^-)$.
\end{enumerate}
\end{notation}

In the present subsection, with $H$ abelian, $H^-$ is just another
name for $T[w]$.  Later on, however, we will use the same notation in
a more general setting, where it must be taken more seriously.

The next lemma is fundamental, and is the only one in which the
precise choice of the involution $w$ is fully exploited. It will be
applied repeatedly. 

\begin{lemma}\label{invertH-}
Under our present hypotheses and notation, 
but without requiring $\hat H$ to be
nonabelian, there is no involution $i\in B$ inverting $H^-$.
\end{lemma}
\begin{proof}
Supposing on the contrary that the involution $i\in B$ inverts $H^-$,
it follows in particular that $H^-$ is abelian.

As $H$ contains no involutions, it also follows that 
$H^-=[i,H^-]\le F\oo(B)$.  Now by the choice of $w$, $\rk(H^-)\ge
\rk(B/C(i))\ge \rk (F\oo(B)/C_{F\oo(B)}(i))$. 
Furthermore, $F(B)$ is a $2^\perp$-group (Lemma \ref{F(B)}).

So as $i$ acts on $F\oo(B)$, the latter 
decomposes definably as a product (of sets) as
$F\oo(B)=C_{F\oo(B)}(i)\cdot F\oo(B)^-$ where $F\oo(B)^-=\{a\in
F\oo(B):a^i=a^{-1}\}$ \cite[Ex.~14, p.~73]{BN}.
As $F\oo(B)$ is
connected, each factor has Morley degree one.
Considering the ranks, it follows that $H^-$ is a generic subset of
$F\oo(B)^-$, and is therefore the unique definable subgroup 
of this rank contained in $F\oo(B)^-$. 
It follows that $C(i)$
normalizes the group $H^-$, 
and in particular some Sylow 2-subgroup $S$ of $B$
normalizes $H^-$, so $N\oo(H^-)\le B$.
By conjugation also $S^w$ normalizes $H^-$.
Hence $S^w\le B$, and this gives a contradiction.
\end{proof}

Next we will give a companion lemma that goes in the opposite
direction: there are involutions in $B$ 
inverting large pieces of $H^-$ (when $H^-$ is abelian). This
is Lemma \ref{invertP} below. We prepare the way with a very general
lemma, which depends on the following fact from the theory of generalized Sylow
subgroups.   

\begin{fact}[{\cite[4.19,4.20,4.22]{Bu-Th}, 
\cite[Theorem 6.7, Cors.~6.8,6.9]{Bu-GCS}}] 
\label{U0r-Sylow}
Let $H$ be a connected solvable group of finite Morley rank, and
$r>0$. Then the Sylow $U_{0,r}$-subgroups $U$ of $H$ are of the following
form
$$U=Q_r\cdot H'_r$$
where $Q$ is a Carter subgroup of $H$, $Q_r=U_{0,r}(Q)$, and
$H'_r=U_{0,r}(H')$. In particular, $U$ is normalized by a Carter
subgroup,
and if $U_{0,r}(H')=1$ then $U$ is contained in a Carter subgroup.
\end{fact}

\begin{lemma}\label{invariance}
Let $H$ be a group of finite Morley rank with $H\oo$ solvable and 
$U_2(H)=1$, 
and let $Q$ be either a Carter subgroup, or a Sylow
$U_{0,r}$-subgroup, of $H\oo$.
Then $N_H(Q)$ contains a Sylow 2-subgroup of $H$.
In particular, every involution of $H$ is conjugate under $H\oo$ to
one which normalizes $Q$.
\end{lemma}
\begin{proof}
We suppose first that $Q$ is a Carter subgroup of $H\oo$.
By Lemma \ref{carter-sylow}, $N_H(Q)$ contains a Sylow $2$-subgroup of
$H\oo$. By the Frattini argument, $N_H(Q)$ covers $H/H\oo$.
If $S$ is a Sylow $2$-subgroup of $N_H(Q)$,
then by Lemma \ref{sylow:basics}, the group $S$ covers a Sylow
$2$-subgroup of $H/H\oo$. 
It follows easily that $S$ is a Sylow $2$-subgroup of $H$.

Now suppose that $Q$ is a Sylow $U_{0,r}$-subgroup of $H\oo$.
Then $Q=Q_1Q_2$ where $Q_1=U_{0,r}([H\oo]')$ and $Q_2=U_{0,r}(C)$ for
some Carter subgroup $C$ of $H\oo$ (Fact \ref{U0r-Sylow}).
Hence $N(C)\le N(Q)$ and the second claim follows.

The final claim is then immediate.
\end{proof}

We need the next lemma, at the moment, under the hypothesis that $H$
is abelian, but it will be applied more generally in the following
section. 

\begin{lemma}\label{invertP}
With our usual hypotheses and notations, but with no additional
assumptions on $\hat H$ or $H$, 
suppose that $H^-$ is abelian.
Let $P=U_{0,r}(H^-)$ for some $r$. 
Then there is an involution $w_P$ in $B$ which inverts $P$.
\end{lemma}
\begin{proof}
We may suppose that $P$ is nontrivial.

Applying Corollary \ref{invariantBorel},
let $\tilde B_0$ be a $w$-invariant Borel subgroup containing
$N\oo(P)$.  
Let $\tilde B_1$ be a Borel subgroup distinct from $B$
with $(B\intersect \tilde B_0)\oo\le \tilde B_1$, chosen so as to maximize
$(B\intersect \tilde B_1)\oo$. Let $\tilde H=(B\intersect \tilde B_1)\oo$.

Suppose $\tilde H$ is {\em abelian}. 
Let $U=U_{0,r}(\tilde H)$. 
Then $N(U)\le N(B)$ by Lemma \ref{QandB} (2).
Hence $U_{0,r}(N_{\tilde B_0}(U))\le U_{0,r}(\tilde H)=U$, and by
Fact \ref{unipotencebasics} (\ref{ub4}) we find that
$U$ is a Sylow $U_{0,r}$-subgroup of
$\tilde B_0$. Accordingly $U$ is also a Sylow $U_{0,r}$-subgroup of
$C\oo(P)$. By Lemma \ref{invariance},
it follows that $U$ is normalized by an involution $w_P$
conjugate to $w$ under the action of $C\oo(P)$.  Hence $w_P$ inverts
$P$ as well.  Since $N(U)\le N(B)$,
the involution $w_P$ normalizes $B$, and hence lies in $B$.

Suppose now that
$\tilde H$ is {\em nonabelian}. 

By Lemma \ref{uniqueness-1} 
we can apply Fact \ref{JL:H'} freely,
which requires also bearing in mind Lemma \ref{9.2Carter'}.

Let $\tilde r'=\rred(\tilde H')$. Suppose first that $r\ne \tilde r'$.
Then by Fact \ref{U0r-Sylow}, $P$ is contained in a Carter subgroup
$Q$ of $\tilde H$.  Now $B$ contains $N\oo(Q)$ by Lemma \ref{QandB}.
It follows that $N_{\tilde B_0}\oo(Q)\le (B\intersect \tilde
B_0)\oo\le (B\intersect \tilde B_1)\oo$, so $N_{\tilde B_0}\oo(Q)\le
N_{\tilde H}(Q)=Q$.  Hence $Q$ is a Carter subgroup of $\tilde B_0$ by
Fact \ref{nearcarter}, and hence also of $C\oo(P)$, as $Q$ is abelian
(Fact \ref{JL:H'}).  So there is an involution $w_P$ conjugate to $w$
under the action of $C\oo(P)$ such that $w_P$ normalizes $Q$.  Since
$w$ inverts $P$, also $w_P$ inverts $P$.  By Lemma \ref{QandB}, $w_P$
normalizes $B$.  Thus in this case we have our claim.

Now suppose that $r=\tilde r'$. Then $U_{0,r}(\tilde H)$ is abelian by
Fact \ref{JL:H'}, and hence is contained in $C\oo(P)\le \tilde B_0$.  On
the other hand, $\tilde H<N_B\oo(F_r(\tilde H))$ by Fact \ref{JL:H'}
(\ref{H':5}),
and thus $N\oo(F_r(\tilde H))$ is contained in $B$, and in no other
Borel subgroup. 
So $N_{\tilde B_0}\oo(F_r(\tilde H))\le (B\intersect \tilde B_0)\oo\le
\tilde H$, and thus
$F_r(\tilde H)=U_{0,r}(\tilde H)$ (Fact \ref{JL:H'})
is a Sylow $U_{0,r}$-subgroup of $B_0$, and hence also
of $C\oo(P)$.  Now by Lemma \ref{invariance} 
it follows that $w$ is conjugate under the action
of $C\oo(P)$ to an involution $w_P$ normalizing $F_r(\tilde H)$.

As $B$ is the only Borel subgroup containing $N\oo(F_r(\tilde H))$,
the involution $w_P$ normalizes $B$, and hence $w_P$ lies in $B$.
So again we have our claim.
\end{proof}

Now we can wrap up the first phase of our analysis.  We will make use
of another two points from the theory of maximal pairs, from
\cite{Bu-JL}.

\begin{fact}[{\cite[9.2 (5d,5b,5c)]{Bu-Th},
\cite[Lemmas 3.12, 3.13]{Bu-JL}}]
\label{JL:Fr(B)}
Let $G$ be a minimal connected simple group of finite Morley rank,
and $(B,B_1)$ a maximal pair of Borel subgroups of $G$. Let
$H=(B\intersect B_1)\oo$. Suppose that $H$ is nonabelian, 
and that $\rred(B)<\rred(B_1)$. 
Then the following hold.
\begin{enumerate}
\item $F_r(B)\le Z(H)$ for $r\ne \rred(H')$.
\item\label{JL:Fr(B)2} $F_{\rred(H')}(B)$ is nonabelian.
\end{enumerate}
In particular, the group $F_r(B)$ is abelian if and only if $r\ne
\rred(H')$. 
\end{fact}

Item $(\ref{JL:Fr(B)2})$ was not noted explicitly in \cite[9.2]{Bu-Th},
but follows readily, and is given in \cite[Theorem 4.5, (4)]{Bu-JL}.

\begin{lemma}
Under our standing hypotheses, 
the Borel subgroup $B_1$ can be chosen to be $w$-invariant.
\end{lemma}
\begin{proof}
In view of  Lemma \ref{Hnonabelian}, we may suppose that 
$$\hbox{$H$ is abelian}$$ 
We will show that in this case 
there is a choice of $B_1$ for which $(B\intersect B_1)$ is abelian,
and thus by Lemma \ref{Hhat-abelian} we may also choose such a $B_1$ 
which is $w$-invariant.
Recall that $H$ is torsion free. 

Suppose toward a contradiction that for all suitable choices of $B_1$, 
$$\hbox{$(B\intersect B_1)\oo$ is nonabelian}$$ 

We first consider any 
$\tilde B_1$ arbitrarily which meets our  basic conditions, and set
$\tilde H=(B\intersect \tilde B_1)\oo$.
As this group is assumed nonabelian,
Fact \ref{JL:Fr(B)} applies in view of Lemma \ref{uniqueness-1}.
In particular, the value of $\rred({\tilde H}')$ is determined by the
structure of $B$, as the value of $r$ for which $F_r(B)$ is
nonabelian.  We will denote this value by $\hat r'$, as usual, but we
emphasize that its value depends only on $B$. 

Now if 
$$\hbox{$U_{0,r}(H^-)=1$ for all $r\ne \hat r'$}$$
then in particular $H^-=U_0(H^-)$, and
by Lemma \ref{invertP} there is an involution $w_P\in B$
inverting $H^-$. This contradicts Lemma \ref{invertH-}.
So, in fact, 
$$\hbox{$U_{0,r}(H^-)\ne 1$ for some $r\ne \hat r'$}$$  

In this case, let 
$P=U_{0,r}(H^-)$ and let $B_0$ be a $w$-invariant Borel subgroup
containing $N\oo(P)$ (Corollary \ref{invariantBorel}). 
Note that $B_0\ne B$.
Let $H_0=(B\intersect B_0)\oo$, and choose $B_1$
distinct from $B$ and containing $H_0$ so that $\hat H=(B\intersect
B_1)\oo$ is maximal.  By our hypothesis, 
$\hat H$ is nonabelian.

By Lemma \ref{invertP}, there is an involution $w_P\in B$ which
inverts $P$. As $P$ is torsion free it is $2$-divisible,
so $P=[w_P,P]\le F(B)$. As $r\ne \hat r'$,
it follows from Fact \ref{unipotencebasics} (\ref{nilstructure})
that $P$ and $F_{\hat r'}(B)$ commute.
Hence $F_{\hat r'}(B)\le C\oo(P)\le  B_0$ and thus
$F_{\hat r'}(B)\le \hat H$,
and this contradicts Facts \ref{JL:H'} (2,7), since we assumed 
that $\hat H$ is nonabelian. 

\end{proof}

At this point, we can take up the analysis afresh. As noted, we will
reuse some of the auxiliary information found along the way
(which has been stated in sufficient generality to allow this), and
from this point on the logic of the argument is completely linear.


%% file: min06.tex

\section{Case II, conclusion}

We recall the notations: $B$ is a standard Borel subgroup with
$N(B)$ strongly embedded in $G$, and $w$ is an involution outside
$N(B)$. $T[w]$ is the set of elements of $B$ inverted by $w$, and by
hypothesis has rank at least $\rk(B/C(i))$ for $i$ any involution of
$B$. $H=(B\intersect B^w)$ (which is torsion free, and in particular
connected). The group $B_1$ is a Borel subgroup 
distinct from $B$, containing $H$, and chosen to maximize $\hat
H=(B\intersect B_1)\oo$. The group $B_1$ is also $w$-invariant.
We fix a Carter subgroup $Q$ of $\hat H$.
Then $Q$ is also a Carter subgroup of $B_1$, and $N(Q)\le N(B)$ by
Lemma \ref{QandB}.
The Pr\"ufer 2-rank is assumed to be at least two.

We have not determined whether or not 
$\hat H$ is abelian, and we will 
frequently have to argue according to cases.
When $\hat H$ is nonabelian, we use the structural information afforded
by \cite{Bu-JL}, specifically Facts \ref{9.2Carter}, \ref{JL:H'},
\ref{JL:Fr(B)}. 
This is justified by Lemma \ref{uniqueness-1}. When $\hat
H$ is abelian we will have to argue directly. In either case we will
arrive at much the same conclusions. 

\subsection{The involution $w_1$}

There is one more essential ingredient in this configuration, as follows.

\begin{lemma}
There is an involution $w_1\in B$ which normalizes $B_1$ and $Q$, and
which is conjugate to $w$ under the action of $B_1$.
\end{lemma}
\begin{proof}
As $Q$ is a Carter subgroup of $B_1$ and $w$ normalizes $B_1$, 
by Lemma \ref{invariance} there is $w_1$ conjugate to $w$
under the action of $B_1$ which normalizes $Q$. Then $w_1$ normalizes
$B_1$, and as $N(Q)\le N(B)$ also $w_1$ normalizes $B$, and hence lies
in $B$.
\end{proof}

We will make use of the following general principle.

\begin{fact}[{\cite[3.18]{Bu-Th}, \cite[3.6]{Bu-SF}}]
\label{centralizer}
Let $H$ be a nilpotent group of finite Morley 
rank and $P$ a group of definable automorphisms of $H$. 
Suppose that
$P$ is a finite $p$-group and
$H$ is a $U_{0,r}$-group with no elements of order $p$.
Then $C_H(P)$ is a $U_{0,r}$-group.
\end{fact}

\begin{lemma}\label{w1onB1}
\ 
\begin{enumerate}
\item $\rred(B_1)>\rred(\hat H)$.
\item The involution $w_1$ inverts $U_0(B_1)$.
\end{enumerate}
\end{lemma}
\begin{proof}
\ 

1. If $\hat H$ is nonabelian, the first claim is given by 
Lemma \ref{9.2Carter'}. 

Suppose now that $\hat H$ is abelian. 
Let $U=U_0(\hat H)$.
Recall that $H$ is torsion free and hence
$\rred(\hat H)\ge \rred(H)>0$.  By Lemma \ref{QandB}, 
$N\oo(U)\le B$. So $N_{B_1}\oo(U)\le \hat H$. It follows that
$\rred(B_1)>\rred(\hat H)$ as Fact \ref{unipotencebasics} (\ref{ub4}) would
apply to $U_0(B_1)$ in the contrary case, implying 
$U=U_0(B_1)$ and $B_1\le N\oo(U)\le B$. 

2. Let $P=U_0(B_1)$.
We consider the action of $w_1$ on $P$.
By Fact \ref{centralizer}, the centralizer in $P$
of $w_1$ is a $U_{0,r_1}$-group with $r_1=\rred(B_1)$. This 
centralizer is contained in $B$, by strong embedding, and since 
$r_1>\rred(\hat H)$, it must be trivial.  Now $P$ is a
$2^\perp$-group,
as otherwise its Sylow $2$-subgroup is central in $B_1$, and thus
$C(i)$ is a Borel subgroup for each involution $i$.
Thus $w_1$ inverts $P$ \cite[p.~78, Ex.~13]{BN}.
\end{proof}

\subsection{$H^-$ and $\hat H^-$}\ 

We now consider more closely the action of $w$ on $H$, and the action
of $w_1$ on $\hat H$. 

\begin{notation}
\ 
\begin{enumerate}
\item Fix an involution $w_1\in B$ normalizing $B_1$ and $Q$.
\item Set $\hat H^-=\<\{a^2:a\in \hat H, a^{w_1}=a^{-1}\}\>$.
\end{enumerate}
\end{notation}

\begin{lemma}\label{hatH-}
The group $\hat H^-$ is an abelian subgroup of $F(B)$,
and contains no involutions. 
\end{lemma}
\begin{proof}
For $a\in \hat H$ with
$a^{w_1}=a^{-1}$, we have
$a^2=[w_1,a]\in F(B)$. Thus $\hat H^-\le F(B)$.

So $\hat H^-\le F(B)\intersect \hat H \le F(\hat H)$,
which is an abelian group, applying Fact \ref{JL:H'} if $\hat H$ is
not itself abelian.  

So $\hat H^-$ is an abelian subgroup of $F(B)$,
inverted by $w_1$. As $F(B)$ contains no involutions, the same applies
to $\hat H^-$.
\end{proof}

It follows in particular that any definable subgroup of
$\hat H^-$ is $2$-divisible, and that
$$\hat H^-=\{a^2:a\in \hat H, a^{w_1}=a^{-1}\};$$
in other words, this set is already a group. 

The structure of $\hat H^-$ is clarified by Lemma 
\ref{U0rhatH-} below. As preparation, we insert an additional
lemma about our maximal pair $(B,B_1)$.

\begin{lemma}\label{U0r-hatH}
For any $r$, and any nontrivial 
Sylow $U_{0,r}$-subgroup $P$ of $\hat H$,
$N\oo(P)\le B$.
\end{lemma}
\begin{proof}
Extend $P$ to a Sylow $U_{0,r}$-subgroup $U$ of $B$.  

We show first:
$$N_B\oo(P)\not \le \hat H\leqno(1)$$

If $P<U$ then the claim is clear by Fact \ref{unipotencebasics}
(\ref{ub4}). 

If $P=U$, then $P$ is normalized by a Carter
subgroup of $B$ (Fact \ref{U0r-Sylow}), and 
a Carter subgroup of $B$ cannot be contained in $\hat H$ (Lemma
\ref{9.2Carter'}). 
This proves $(1)$.  

Now if $\hat H$ is abelian then $N_B\oo(P)>\hat H$ and the lemma
follows by the maximality of $\hat H$. 
So we suppose
$$\hbox{$\hat H$ is nonabelian}\leqno(2)$$
In this case Facts \ref{9.2Carter}, \ref{JL:H'}, and
\ref{JL:Fr(B)} apply.
We set $\hat r'=\rred(\hat H')$. 

If $r\ne \hat r'$ then $F_r(B)\le P$ (Fact \ref{JL:Fr(B)}).  
We know that $U$ is normalized
by a Carter subgroup of $B$, and in particular (Lemma \ref{carter-sylow})
by a Sylow $2$-subgroup
$S$ of $B$.  In particular $\Omega_1(S)$ normalizes $U$.  For $i$ an
involution in $S$, we have $[i,U]\le F(B)\intersect U$, that is $i$
acts trivially on $U/(F(B)\intersect U)$ and $F(B)\intersect U$ is a 
$2^\perp$-group.  By \cite[3.2]{Bu-SF}, 
the centralizer $C_U(i)$ then covers the quotient, that is
we have $U=C_U(i)(F(B)\intersect U)$, and so
by Fact \ref{U0r-generation}, 
we have $U=C_U(i)F_r(B)$. Hence $[i,P]\le [i,U]\le F_r(B)\le P$ and
$i$ normalizes $P$. Thus $\Omega_1(S)$ normalizes $P$, and as the
Pr\"ufer $2$-rank is at least two it follows that $N\oo(P)\le B$
(Lemma \ref{blackhole}),  as claimed.

If $r=\hat r'$ then $P=U_{0,r}(\hat H)=F_r(\hat H)$ by Fact
\ref{JL:H'} (2), so $\hat H<N_B\oo(P)$ by Fact \ref{JL:H'} (7) and
Fact \ref{unipotencebasics} (5),
and our claim follows again by maximality of $\hat H$.
\end{proof}

Now we can control the structure of $\hat H^-$. We will make use of
the following very general result.

\begin{lemma}\label{inversion-commuting}
Let $G$ be a group, and $H,K\le G$ subgroups with $K$ normalizing
$H$. Let $t\in G$ act on $H$ and $K$, inverting both groups. If $K$ is
$2$-divisible, then $[K,H]=1$.
\end{lemma}
\begin{proof}
For $h\in H$ and $k\in K$ we have
$$(h^{-1})^k=(h^k)^{-1}=(h^k)^t=(h^t)^{k^t}=(h^{-1})^{k^{-1}}$$ 
and thus
$(h^{-1})^{k^2}=h^{-1}$.
\end{proof}

\begin{lemma}\label{U0rhatH-}
If $U_{0,r}(\hat H')=1$ then $U_{0,r}(\hat H^-)=1$.
\end{lemma}

\begin{proof}
Let $P=U_{0,r}(\hat H^-)$.  Suppose $P>1$.

The involution $w_1$ inverts $U_0(B_1)$ and $P$; the
latter is $2$-divisible.  It follows by Lemma \ref{inversion-commuting} 
that $P$ centralizes
$U_0(B_1)$. 

Suppose first that $\hat H$ is noncommutative. 
Then by Fact \ref{JL:Fr(B)},
since $P\le F_r(B)\le Z(\hat H)$, 
the group $P$ commutes with $F(B)$
in view of the structure theory of Fact \ref{unipotencebasics}
(\ref{nilstructure}). 
So $C\oo(P)$ contains $F(B)$, $U_0(B_1)$, and $\hat H$,
and hence either $F(B)\le B_1$ or $U_0(B_1)\le B$. But
$\rred(B_1)>\rred(B)$, so the second possibility is excluded,
and the first possibility is excluded by Fact \ref{JL:H'} (2,7).

So we may now suppose that 
$$\hbox{$\hat H$ is commutative.}$$

We claim: 
$$\hbox{\em For any $s$, either $F_s(B_1)\le \hat H$, or $U_{0,s}(\hat
H)=1$.}\leqno(*)$$ 
Suppose that $F_s(B_1)\not \le \hat H$, and let $U$
be a Sylow $U_{0,s}$-subgroup of $\hat H$. Then the group $F_s(B_1)U$
is nilpotent by Fact \ref{UbyU}.  So $N_{B_1}\oo(U)\not \le \hat H$
and by Lemma \ref{U0r-hatH} we have $U=1$. Our claim $(*)$ follows.

Now $\hat H$ is a Carter subgroup of $B_1$, by Lemma 
\ref{QandB} (1).
So $B_1=F(B_1)\hat H$
\cite{Wa-NCC}.
Furthermore, as $\hat H$ is commutative, $H^-$
centralizes any factor $F_s(B_1)$ which lies in $\hat H$.

Consider a factor $F_s(B_1)$, where $U_{0,s}(\hat H)=1$.  The
centralizer in $F_s(B_1)$ of $w_1$ is a $U_{0,s}$-group, by Fact
\ref{centralizer}, and is contained in $\hat H$, hence is trivial.  
As $F(B_1)$ is a $2^\perp$-group, it follows that
$w_1$ inverts $F_s(B_1)$ \cite[Ex.~14, p.~73]{BN}. 
As $w_1$ and $w$ are conjugate under the
action of $B_1$, also $w$ inverts $F_s(B_1)$.  But $w$ also inverts
$H^-$, and as $H^-\le H$ is $2$-divisible it follows that $H^-$
commutes with $F_s(B_1)$ by Lemma \ref{inversion-commuting}.

So $H^-$
commutes with every factor $F_s(B_1)$ for $s\ge 1$.

Furthermore, since $H^-$ is torsion free,
it commutes with every $U_p(B_1)$ for $p$ a prime,
as a consequence of the main result of \cite{Wa-FFM}; this is given
in \cite[3.13]{AlCh-LSE} for actions on abelian unipotent
$p$-groups, and the general case has the same proof.

Since the divisible part of the torsion subgroup of $F(B_1)$ is
central in $B_1$, 
it follows that $H^-$ centralizes $F\oo(B_1)$. 
Hence $H^-$ centralizes $F\oo(B_1)\hat H=B_1$, and $H^-\le Z(B_1)$.

Now $w$ and $w_1$ are conjugate under the action of $B_1$, and
therefore $w_1$ inverts $H^-$. 
This contradicts Lemma \ref{invertH-}.
\end{proof}

\subsection{Invariants attached to $w$}\ 

Now Lemma \ref{U0rhatH-} produces a peculiar situation. 
If, for example, $\hat H$ is abelian, then it follows that
$U_{0,r}(\hat H^-)=1$ for all $r$. On the other hand, this is certainly not
the case for $H^-$, which is torsion free. Furthermore $w$ and $w_1$
are conjugate under the action of $B_1$. Of course, this conjugation
need not preserve $\hat H$, or carry $H^-$ to $\hat H^-$, so 
we are still short of a contradiction. However, it is possible to
attach certain invariants to involutions acting on solvable
groups, which will be preserved by conjugation, and in this way arrive
at a contradiction by comparing the actions of $w$ and $w_1$.

This is based on the following considerations.

\begin{lemma}\label{invariantsylow}
Let $H$ be a solvable group of finite Morley rank with $H\oo$
$2$-divisible, and $w$ an involution in $H$. Fix $r\ge 1$. Then the
following hold.
\begin{enumerate}
\item Any $w$-invariant nilpotent $U_{0,r}$-subgroup 
$P$ of $H$ is contained in a $w$-invariant Sylow 
$U_{0,r}$-subgroup of $H$.
\item 
If a Sylow $2$-subgroup of $H$ contains a unique involution,
then any two $w$-invariant Sylow $U_{0,r}$-subgroups of $H$ are
conjugate under the action of $C_H(w)$. 
\end{enumerate}
\end{lemma}
\begin{proof}
1. We may suppose that $P$ is a maximal $w$-invariant
$U_{0,r}$-subgroup of $H$. We claim that $P$ is a Sylow
$U_{0,r}$-subgroup of $H$.  It suffices to show that $P$ is a Sylow
$U_{0,r}$-subgroup of $N_H\oo(P)$, or in other words we may suppose
that $H$ normalizes $P$.  In this case, replacing $H$ by $H/P$, we may
suppose that $P=1$.
In this case, our claim reduces to Lemma \ref{invariance}.

2. Let $P$ and $P^g$ (with $g\in H$) be two $w$-invariant
Sylow $U_{0,r}$-subgroups of $H$. 

Then $w$ and $w^g$ are in $N_H(P^g)$, and the latter group has a
unique involution in each Sylow $2$-subgroup. Hence
$w$ and $w^g$ are conjugate in $N_H(P^g)$. If $w^{gh}=w$ with $h\in
N_H(P^g)$, then $P^g=P^{gh}$ and $gh\in C_H(w)$.
\end{proof}

\begin{definition}
Let $w$ be an involution acting definably on a solvable group $H$ of
finite Morley rank, and let $r\ge 1$. 
Then we define
$\iota_r(w,H)$ as the maximal rank of a $U_{0,r}$-subgroup of $H$
inverted by $w$.
\end{definition}

Note that one would normally require the $U_{0,r}$-subgroups involved
to be nilpotent, but here they are in any case abelian.

What makes this invariant manageable is the following.

\begin{lemma}
Let $w$ be an involution acting definably on a solvable group $H$ of
finite Morley rank with $U_2(H)=1$,  and let $r\ge 1$. 
Suppose that a Sylow $2$-subgroup of $H\<w\>$ contains a unique
involution.  
Let $U$ be a $w$-invariant Sylow
$U_{0,r}$-subgroup of $H$. Then $\iota_r(w,H)=\iota_r(w,U)$. 
\end{lemma}
\begin{proof}
Let $P$ be a $U_{0,r}$-subgroup of $H$ inverted by $w$, and of maximal
rank. Extend $P$ to a Sylow $U_{0,r}$-subgroup $Q$ of $H$. 
Then $Q$ and $U$ are conjugate under the action of $C(w)$, and our
claim follows.
\end{proof}

Let us now specialize this to the case at hand.

\begin{lemma}\label{iota}
With our usual hypotheses and notation, 
let $U_1$ be a $w_1$-invariant Sylow $U_{0,r}$-subgroup of $B_1$.
Then $\iota_r(w_1,U_1)=\iota_r(w,B_1)$.
\end{lemma}
\begin{proof}
Let $\hat B_1=B_1\<w\>=B_1\<w_1\>$. 
We show first that 
$$\hbox{A Sylow $2$-subgroup of $\hat B_1$ contains a unique
  involution.}\leqno(*)$$ 

If $B_1$ is a $2^\perp$-group, then the Sylow subgroups of $\hat
B_1$ are cyclic of order two. 

Suppose $B_1$ contains an involution. Then $B_1$ meets a conjugate
$M_1$ of $M$. 

If $\hat B_1$ is contained in $M_1$, then $B_1$ is conjugate
to $B$. But this case may be ruled out as follows.
By Lemma \ref{w1onB1}, $w_1$ inverts
$U_0(B_1)$. It follows that some involution of $B$ inverts $U_0(B)$,
and hence all involutions of $B$ invert $U_0(B)$, which is impossible.

So $\hat B_1$ meets $M_1$ in a proper subgroup of $\hat B_1$, 
which is therefore
strongly embedded in $\hat B_1$, and in particular all involutions of
$\hat B_1$ are conjugate, and lie in $B_1$. So if $\hat B_1$ contains
an elementary abelian $2$-group of order $4$, then the same applies
to $B_1$, and hence by Lemma \ref{blackhole}, $B_1$ is conjugate to
$B$, which we have just ruled out. So $(*)$ holds in all cases.

Now the general theory applies to $w$ and $w_1$ in $\hat B_1$,
and as they are conjugate we find
$$\iota_r(w, B_1)=\iota_r(w,\hat B_1)=\iota_r(w_1,\hat B_1)=\iota_r(w_1,U_1)$$
\end{proof}

\subsection{Structure of $H^-$}\ 

Now we take up the structure of $H^-$, about which we
know very little at this point; recall that this group is the
definable closure of the group generated by $T[w]$.

\begin{lemma}\label{U0r(H-)}
If $U_{0,r}(\hat H')=1$ then $U_{0,r}(H^-)=1$.
\end{lemma}
\begin{proof}
Let $P$ be a $w$-invariant Sylow $U_{0,r}$-subgroup of $H^-$, and 
suppose that $P$ is nontrivial. We claim:
$$\hbox{$w$ inverts $P$}\leqno(1)$$

If $H^-$ is abelian, then $w$ inverts $H^-$ and this is clear.

If $H^-$ is nonabelian, then also $\hat H$ is nonabelian, and as usual
we set $\hat r'=\rred({\hat H}')$. By assumption $r\ne \hat r'$,
so $U_{0,r}({\hat H}')=1$ by Fact \ref{JL:H'} (1).
Let $\overline {H^-}=H^-/(H^-)'$.  Then
$\overline{H^-}$ is an abelian group inverted by $w$, and
$U_{0,r}([H^-]')=1$.

Let $P_0=C_P(w)$. Then $P_0$ is a $U_{0,r}$-group by Fact
\ref{centralizer}. So the image $\bar P_0$ of $P_0$ in
$\overline{H^-}$ is a $U_{0,r}$-group which is both centralized and
inverted by $w$. Hence $\bar P_0=1$ and  $P_0\le (H^-)'$. 
Since $P_0$ is a $U_{0,r}$-group,
it follows that $P_0=1$. So $w$ inverts $P$, and $(1)$ holds
in either case.

Let $P_1$ be a $w_1$-invariant Sylow $U_{0,r}$-subgroup of $\hat H$.
By Lemma \ref{U0r-hatH}, $N\oo(P_1)$ is contained in $B$ and hence
$P_1$ is a Sylow $U_{0,r}$-subgroup of $B_1$. 
By Lemma \ref{iota}, we have $\iota_r(w_1,P_1)=\iota_r(w,B_1)>0$. 
This contradicts Lemma \ref{U0rhatH-}.
\end{proof}

Now the analysis of $H^-$ produces a contradiction.

\begin{lemma}\label{nonabelian}
This configuration is inconsistent.
\end{lemma}
\begin{proof}
If $\hat H$ is abelian, then the preceding lemma shows that
$U_{0,s}(H^-)=1$ for all $s$, and as $H^-$ is torsion free and
nontrivial this is a contradiction.

So $\hat H$ is nonabelian, and 
we may apply Fact \ref{JL:H'}.
As usual we set $\hat r'=\rred(\hat H')$.
Then by Lemma \ref{U0r(H-)},
$U_{0,r}(H^-)=1$ for $r\ne \hat r'$. 
As $H^-\le H$ is torsion free, it follows that 
$H^-=U_{0,\hat r'}(H^-)\le U_{0,\hat r'}(\hat H)=F_{\hat r'}(\hat H)$,
which is abelian. 

As $H^-$ is abelian, we can apply
Lemma \ref{invertP}, and there is an involution in $B$ which inverts
$H^-$, which contradicts Lemma \ref{invertH-}.
\end{proof}

With this contradiction, the proof of Theorem \ref{main} in the second case is
finally complete.